\input amstex
\documentstyle{amsppt}
\magnification=\magstep1
\vsize =21 true cm
\hsize =16 true cm
\loadmsbm
\topmatter

\centerline{\bf Dedekind $\eta$-function, Hauptmodul and invariant
                theory}
\author{\smc Lei Yang}\endauthor
\endtopmatter
\document

\centerline{\bf Abstract}

\vskip 0.5 cm

  We solve a long-standing open problem with its own long history
dating back to the celebrated works of Klein and Ramanujan. This
problem concerns the invariant decomposition formulas of the
Hauptmodul for $\Gamma_0(p)$ under the action of finite simple
groups $PSL(2, p)$ with $p=5, 7, 13$. The cases of $p=5$ and $7$
were solved by Klein and Ramanujan. Little was known about this
problem for $p=13$. Using our invariant theory for $PSL(2, 13)$,
we solve this problem. This leads to a new expression of the
classical elliptic modular function of Klein: $j$-function in
terms of theta constants associated with $\Gamma(13)$. Moreover,
we find an exotic modular equation, i.e., it has the same form as
Ramanujan's modular equation of degree $13$, but with different
kinds of modular parametrizations, which gives the geometry of the
classical modular curve $X(13)$.

\vskip 0.5 cm

\centerline{\bf Contents}
$$\aligned
 &\text{1. Introduction}\\
 &\text{2. Six-dimensional representations of $PSL(2, 13)$ and transformation formulas}\\
 &\text{\quad for theta constants}\\
 &\text{3. Seven-dimensional representations of $PSL(2, 13)$, exotic modular equation and}\\
 &\text{\quad geometry of modular curve $X(13)$}\\
 &\text{4. Fourteen-dimensional representations of $PSL(2, 13)$ and invariant decomposition}\\
 &\text{\quad formula}\\
 \endaligned$$

\centerline{\bf 1. Introduction}

\vskip 0.5 cm

  In many applications of elliptic modular functions to number theory
the Dedekind eta function plays a central role. It is defined in
the upper-half plane ${\Bbb H}=\{ z \in {\Bbb C}: \text{Im}(z)>0
\}$ by
$$\eta(z):=q^{\frac{1}{24}} \prod_{n=1}^{\infty} (1-q^n), \quad q=e^{2 \pi i z}.$$
The Dedekind eta function is closely related to the partition
function. A partition of a positive integer $n$ is any
non-increasing sequence of positive integers whose sum is $n$. Let
$p(n)$ denote the number of partitions of $n$. The partition
function $p(n)$ has the well-known generating function
$$\sum_{n=0}^{\infty} p(n) q^n=\prod_{n=1}^{\infty} \frac{1}{1-q^n}
 =\frac{q^{\frac{1}{24}}}{\eta(z)}.$$

  In his ground-breaking works \cite{R1}, \cite{R2} and \cite{R3}, Ramanujan
found the following famous Ramanujan partition congruences:
$$\aligned
  p(5n+4) \equiv 0 \quad &(\text{mod $5$}),\\
  p(7n+5) \equiv 0 \quad &(\text{mod $7$}),\\
 p(11n+6) \equiv 0 \quad &(\text{mod $11$}).
\endaligned$$
Ramanujan's proofs of the congruences modulo $5$ and $7$ are quite
ingenious but are not terribly difficult:
$$\aligned
  q^{\frac{19}{24}} \sum_{n=0}^{\infty} p(5n+4) q^n
 &=5 \frac{\eta(5z)^5}{\eta(z)^6},\\
  q^{\frac{17}{24}} \sum_{n=0}^{\infty} p(7n+5) q^n
 &=7 \frac{\eta(7z)^3}{\eta(z)^4}+49 \frac{\eta(7z)^7}{\eta(z)^8},
\endaligned\eqno{(1.1)}$$
while the proof of the congruence modulo $11$ is much harder. One
of the significance about Ramanujan partition congruences comes
from the problem of classifying congruences of the form $p(ln+b)
\equiv 0$ (mod $l$), where $l$ is a prime. It was proved in
\cite{AB} that there are only those three originally observed by
Ramanujan, when $(l, b)=(5, 4)$, $(7, 5)$, or $(11, 6)$.

  In order to provide a combinatorial explanation of these
congruences, Dyson in \cite{D} defined the rank of a partition to
be its largest part minus the number of its parts. He conjectured
that the partitions of $5n+4$ (resp. $7n+5$) are divided into $5$
(resp. $7$) groups of equal size when sorted by their ranks modulo
$5$ (resp. $7$), thereby providing a combinatorial explanation for
the congruences mod $5$ and $7$. This conjecture was proved by
Atkin and Swinnerton-Dyer (see \cite{ASD}). However, for the third
congruence, the corresponding criterion failed, and so Dyson
conjectured the existence of a statistic, which he called the
crank to combinatorially explain the congruence $p(11n+6) \equiv
0$ (mod $11$). The crank of a partition was found by Andrews and
Garvan (see \cite{AG}). The crank divides the partition into
equinumerous congruence classes modulo $5$, $7$ and $11$ for the
three congruences, respectively. Further work of Garvan, Kim and
Stanton (see \cite{GKS}) produces, for the congruences with moduli
$5$, $7$, $11$ and $25$, combinatorial interpretations which are
rooted in the modular representation theory of the symmetric
group. In fact, in his lost notebook \cite{R}, Ramanujan had
recorded the generating functions for both the rank and the crank.
Moreover, he actually recorded an observation that is equivalent
to Dyson's assertion about the congruences for the ranks of the
partitions of $5n+4$ (see \cite{G1} and \cite{BCCL}). Recently,
Bringmann and Ono (see \cite{BrO}) showed that the rank generating
functions are related to harmonic weak Maass forms. Using
Shimura's theory of half-integer weight modular forms and Serre's
theory of $p$-adic modular forms, they obtained as an application
infinite families of congruences for the rank. These generalize
previous work of Ono (see [O]) on the partition function. In
\cite{Ma}, Mahlburg developed an analogous theory of congruences
for crank generating functions, which are shown to possess the
same sort of arithmetic properties.

  In \cite{Z} Zuckerman found an identity in the spirit of
(1.1) (see also \cite{R} and \cite{BO}):
$$\aligned
  q^{\frac{11}{24}} \sum_{n=0}^{\infty} p(13n+6) q^n
=&11 \frac{\eta(13z)}{\eta(z)^2}+36 \cdot 13
  \frac{\eta(13z)^3}{\eta(z)^4}
 +38 \cdot 13^2 \frac{\eta(13z)^5}{\eta(z)^6}+20 \cdot 13^3 \frac{\eta(13z)^7}{\eta(z)^8}\\
 &+6 \cdot 13^4 \frac{\eta(13z)^9}{\eta(z)^{10}}+13^5 \frac{\eta(13z)^{11}}{\eta(z)^{12}}
  +13^5 \frac{\eta(13z)^{13}}{\eta(z)^{14}}.
\endaligned\eqno{(1.2)}$$
Rademacher (see \cite{Ra}) pointed out that (1.1) can be rewritten
as
$$\aligned
  \sum_{\lambda=0}^{4} \eta(5z) \eta \left(\frac{z+24 \lambda}{5}\right)^{-1}
 &=5^2 \left(\frac{\eta(5z)}{\eta(z)}\right)^6,\\
  \sum_{\lambda=0}^{6} \eta(7z) \eta \left(\frac{z+24 \lambda}{7}\right)^{-1}
 &=7^2 \left(\frac{\eta(7z)}{\eta(z)}\right)^4+7^3
 \left(\frac{\eta(7z)}{\eta(z)}\right)^8.
\endaligned\eqno{(1.3)}$$
In fact, (1.2) can also be rewritten as
$$\aligned
 &\sum_{\lambda=0}^{12} \eta(13z) \eta \left(\frac{z+24 \lambda}{13}\right)^{-1}
 =11 \cdot 13 \left(\frac{\eta(13z)}{\eta(z)}\right)^2+\\
 &+36 \cdot 13^2 \left(\frac{\eta(13z)}{\eta(z)}\right)^4
  +38 \cdot 13^3 \left(\frac{\eta(13z)}{\eta(z)}\right)^6
  +20 \cdot 13^4 \left(\frac{\eta(13z)}{\eta(z)}\right)^8+\\
 &+6 \cdot 13^5 \left(\frac{\eta(13z)}{\eta(z)}\right)^{10}
  +13^6 \left(\frac{\eta(13z)}{\eta(z)}\right)^{12}
  +13^6 \left(\frac{\eta(13z)}{\eta(z)}\right)^{14}.
\endaligned\eqno{(1.4)}$$

  Let $X=X_0(p)$ be the compactification of ${\Bbb H}/\Gamma_0(p)$
where $p$ is prime. The complex function field of $X$ consists of
the modular functions $f(z)$ for $\Gamma_0(p)$ which are
meromorphic on the extended upper half-plane. A function $f$ lies
in the rational function field ${\Bbb Q}(X)$ if and only if the
Fourier coefficients in its expansion at $\infty$: $f(z)=\sum a_n
q^n$ are all rational numbers. The field ${\Bbb Q}(X)$ is known to
be generated over ${\Bbb Q}$ by the classical $j$-functions
$$\left\{\aligned
  j &=j(z)=q^{-1}+744+196884 q+\cdots,\\
  j_p &=j \left(\frac{-1}{pz}\right)=j(pz)=q^{-p}+744+\cdots.
\endaligned\right.$$
A further element in the function field ${\Bbb Q}(X)={\Bbb Q}(j,
j_p)$ is the modular unit $u=\Delta(z)/\Delta(pz)$ with divisor
$(p-1) \{ (0)-(\infty) \}$ where $\Delta(z)$ is the discriminant.
If $m=\text{gcd}(p-1, 12)$, then an $m$th root of $u$ lies in
${\Bbb Q}(X)$. This function has the Fourier expansion
$$t=\root m \of{u}=q^{(1-p)/m} \prod_{n \geq 1} \left(\frac{1-q^n}
    {1-q^{np}}\right)^{24/m}=\left(\frac{\eta(z)}{\eta(pz)}\right)^{24/m}.$$
When $p-1$ divides $12$, so $m=p-1$, the function $t$ is a
Hauptmodul for the curve $X$ which has genus zero (see \cite{Gr}).
It is well-known that the genus of the modular curve $X$ for prime
$p$ is zero if and only if $p=2, 3, 5, 7, 13$. In his paper
\cite{K1}, Klein studied the modular equations of orders $2$, $3$,
$5$, $7$, $13$ with degrees $3$, $4$, $6$, $8$, $14$,
respectively. They are the so-called Hauptmoduln (principal
moduli) which means that there is a uniformising function
$J_{\Gamma}$ that is a modular function for some congruence
subgroup $\Gamma$, and all other modular functions for $\Gamma$
can be written as a rational function in it. Hence, (1.3) and
(1.4) are intimately related to the Hauptmoduln for modular curves
$X_0(5)$, $X_0(7)$ and $X_0(13)$.

  The modular curves $X_0(5)$ and $X_0(7)$ were studied by Klein
in his pioneering work (see \cite{K}, \cite{K1}, \cite{K2},
\cite{KF1} and \cite{KF2}) and by Ramanujan (see \cite{R5},
\cite{R6} and \cite{R}).
  It is well-known that the celebrated Rogers-Ramanujan identities
$$G(q):=\sum_{n=0}^{\infty} \frac{q^{n^2}}{(1-q) \cdots (1-q^n)}
       =\prod_{n=0}^{\infty} \frac{1}{(1-q^{5n+1})(1-q^{5n+4})}$$
and
$$H(q):=\sum_{n=0}^{\infty} \frac{q^{n^2+n}}{(1-q) \cdots (1-q^n)}
       =\prod_{n=0}^{\infty} \frac{1}{(1-q^{5n+2})(1-q^{5n+3})}$$
are intimately associated with the Rogers-Ramanujan continued
fraction
$$R(q):=\cfrac q^{\frac{1}{5}}\\
      1+\cfrac q\\
      1+\cfrac q^2\\
      1+\cfrac q^3\\
      1+\cdots \endcfrac,$$
namely, they satisfy $R(q)=q^{\frac{1}{5}} \frac{H(q)}{G(q)}$. In
\cite{R4}, Ramanujan found an algebraic relation between $G(q)$
and $H(q)$:
$$G^{11}(q) H(q)-q^2 G(q) H^{11}(q)=1+11 q G^6(q) H^6(q),$$
which is equivalent to one of the most important formulas for
$R(q)$ (see also \cite{AB1}):
$$\frac{1}{R^5(q)}-11-R^5(q)=\left(\frac{\eta(z)}{\eta(5z)}\right)^6.\eqno{(1.5)}$$

  In his celebrated work on elliptic modular functions (see
\cite{KF1}, p.640, \cite{KF2}, p.73 and \cite{Du}), Klein showed
that $R(q)=\frac{a(z)}{b(z)}$ where
$$a(z)=e^{-\frac{3 \pi i}{10}} \theta \left[\matrix
       \frac{3}{5}\\ 1 \endmatrix\right](0, 5z), \quad
  b(z)=e^{-\frac{\pi i}{10}} \theta \left[\matrix
       \frac{1}{5}\\ 1 \endmatrix\right](0, 5z).$$
Hence, (1.5) is equivalent to the following formula
$$\left(\frac{\eta(z)}{\eta(5z)}\right)^6=-\frac{f(a(z), b(z))}{a(z)^6 b(z)^6},\eqno{(1.6)}$$
where
$$f(z_1, z_2)=z_1 z_2 (z_1^{10}+11 z_1^5 z_2^5-z_2^{10})$$
is an invariant of degree $12$ associated to the icosahedron,
i.e., it is invariant under the action of the simple group $PSL(2,
5)$, and $z_1^2 z_2^2$ is invariant under the action of the image
of a Borel subgroup of $PSL(2, 5)$ (i.e. a maximal subgroup of
order $10$ of $PSL(2, 5)$) (see \cite{K} for more details). We
call (1.6) the invariant decomposition formula for the
icosahedron.

  In his work on elliptic modular functions (see \cite{KF1}, p.746),
Klein also obtained the invariant decomposition formula for the
simple group $PSL(2, 7)$ of order $168$:
$$\left(\frac{\eta(z)}{\eta(7z)}\right)^4=\frac{\Phi_6(a(z), b(z), c(z))}{a(z)^2 b(z)^2 c(z)^2},\eqno{(1.7)}$$
where
$$a(z)=-e^{-\frac{5 \pi i}{14}} \theta \left[\matrix
       \frac{5}{7}\\ 1 \endmatrix\right](0, 7z), \quad
  b(z)=e^{-\frac{3 \pi i}{14}} \theta \left[\matrix
       \frac{3}{7}\\ 1 \endmatrix\right](0, 7z), \quad
  c(z)=e^{-\frac{\pi i}{14}} \theta \left[\matrix
       \frac{1}{7}\\ 1 \endmatrix\right](0, 7z),$$
and
$$\Phi_6(x, y, z)=xy^5+yz^5+zx^5-5 x^2 y^2 z^2$$
is an invariant of degree $6$ associated to $PSL(2, 7)$, and $xyz$
is invariant under the action of the image of a Borel subgroup of
$PSL(2, 7)$ (i.e. a maximal subgroup of order $21$ of $PSL(2, 7)$)
(see \cite{K2} for more details). Independently, Ramanujan also
gave the same formula (see \cite{R6}, p.300, \cite{B2}, p.174,
\cite{El} and \cite{La} for more details).

  The invariant decomposition formula for the simple group $PSL(2, 13)$ of
order $1092$ is much harder and the modular curve $X_0(13)$ is
much more complicated. In his monograph (see \cite{KF2}, p.73),
Klein showed that $\left(\eta(z)/\eta(13 z)\right)^2$ is a
Hauptmodul for $\Gamma_0(13)$. Ramanujan also studied this problem
(see \cite{R5}, p.326, \cite{R6}, p.244 and \cite{B1}, p.372).
However, neither of them could obtain the invariant decomposition
formula in the spirit of (1.6) and (1.7) (see the end of section 4
for more details).

  In the present paper, we establish the invariant theory for $PSL(2, 13)$.
Combining with theta constants for $\Gamma(13)$, we obtain the
invariant decomposition formula for $PSL(2, 13)$. Let
$$\left\{\aligned
  a_1(z) &:=e^{-\frac{11 \pi i}{26}} \theta \left[\matrix \frac{11}{13}\\ 1 \endmatrix\right](0, 13z)
           =q^{\frac{121}{104}} \sum_{n \in {\Bbb Z}} (-1)^n q^{\frac{1}{2}(13n^2+11n)},\\
  a_2(z) &:=e^{-\frac{7 \pi i}{26}} \theta \left[\matrix \frac{7}{13}\\ 1 \endmatrix\right](0, 13z)
           =q^{\frac{49}{104}} \sum_{n \in {\Bbb Z}} (-1)^n q^{\frac{1}{2}(13n^2+7n)},\\
  a_3(z) &:=e^{-\frac{5 \pi i}{26}} \theta \left[\matrix \frac{5}{13}\\ 1 \endmatrix\right](0, 13z)
           =q^{\frac{25}{104}} \sum_{n \in {\Bbb Z}} (-1)^n q^{\frac{1}{2}(13n^2+5n)},\\
  a_4(z) &:=-e^{-\frac{3 \pi i}{26}} \theta \left[\matrix \frac{3}{13}\\ 1 \endmatrix\right](0, 13z)
           =-q^{\frac{9}{104}} \sum_{n \in {\Bbb Z}} (-1)^n q^{\frac{1}{2}(13n^2+3n)},\\
  a_5(z) &:=e^{-\frac{9 \pi i}{26}} \theta \left[\matrix \frac{9}{13}\\ 1 \endmatrix\right](0, 13z)
           =q^{\frac{81}{104}} \sum_{n \in {\Bbb Z}} (-1)^n q^{\frac{1}{2}(13n^2+9n)},\\
  a_6(z) &:=e^{-\frac{\pi i}{26}} \theta \left[\matrix \frac{1}{13}\\ 1 \endmatrix\right](0, 13z)
           =q^{\frac{1}{104}} \sum_{n \in {\Bbb Z}} (-1)^n q^{\frac{1}{2}(13n^2+n)}
\endaligned\right.\eqno{(1.8)}$$
be the theta constants for $\Gamma(13)$. Set
$$\Phi_{12}(z_1, z_2, z_3, z_4, z_5, z_6)
 :=-\frac{1}{26}(7 \cdot 13^2 {\Bbb G}_0^2+{\Bbb G}_1 {\Bbb G}_{12}
   +{\Bbb G}_2 {\Bbb G}_{11}+\cdots+{\Bbb G}_6 {\Bbb G}_7),\eqno{(1.9)}$$
where
$$\left\{\aligned
  {\Bbb G}_0 &={\Bbb D}_0^2+{\Bbb D}_{\infty}^2,\\
  {\Bbb G}_1 &=-{\Bbb D}_7^2+2 {\Bbb D}_0 {\Bbb D}_1+10 {\Bbb D}_{\infty} {\Bbb D}_1
               +2 {\Bbb D}_2 {\Bbb D}_{12}-2 {\Bbb D}_3 {\Bbb D}_{11}-4 {\Bbb D}_4 {\Bbb D}_{10}
               -2 {\Bbb D}_9 {\Bbb D}_5,\\
  {\Bbb G}_2 &=-2 {\Bbb D}_1^2-4 {\Bbb D}_0 {\Bbb D}_2+6 {\Bbb D}_{\infty} {\Bbb D}_2
               -2 {\Bbb D}_4 {\Bbb D}_{11}+2 {\Bbb D}_5 {\Bbb D}_{10}-2 {\Bbb D}_6 {\Bbb D}_9
               -2 {\Bbb D}_7 {\Bbb D}_8,\\
  {\Bbb G}_3 &=-{\Bbb D}_8^2+2 {\Bbb D}_0 {\Bbb D}_3+10 {\Bbb D}_{\infty} {\Bbb D}_3
               +2 {\Bbb D}_6 {\Bbb D}_{10}-2 {\Bbb D}_9 {\Bbb D}_7-4 {\Bbb D}_{12} {\Bbb D}_4
               -2 {\Bbb D}_1 {\Bbb D}_2,\\
  {\Bbb G}_4 &=-{\Bbb D}_2^2+10 {\Bbb D}_0 {\Bbb D}_4-2 {\Bbb D}_{\infty} {\Bbb D}_4
               +2 {\Bbb D}_5 {\Bbb D}_{12}-2 {\Bbb D}_9 {\Bbb D}_8-4 {\Bbb D}_1 {\Bbb D}_3
               -2 {\Bbb D}_{10} {\Bbb D}_7,\\
  {\Bbb G}_5 &=-2 {\Bbb D}_9^2-4 {\Bbb D}_0 {\Bbb D}_5+6 {\Bbb D}_{\infty} {\Bbb D}_5
               -2 {\Bbb D}_{10} {\Bbb D}_8+2 {\Bbb D}_6 {\Bbb D}_{12}-2 {\Bbb D}_2 {\Bbb D}_3
               -2 {\Bbb D}_{11} {\Bbb D}_7,\\
  {\Bbb G}_6 &=-2 {\Bbb D}_3^2-4 {\Bbb D}_0 {\Bbb D}_6+6 {\Bbb D}_{\infty} {\Bbb D}_6
               -2 {\Bbb D}_{12} {\Bbb D}_7+2 {\Bbb D}_2 {\Bbb D}_4-2 {\Bbb D}_5 {\Bbb D}_1
               -2 {\Bbb D}_8 {\Bbb D}_{11},\\
  {\Bbb G}_7 &=-2 {\Bbb D}_{10}^2+6 {\Bbb D}_0 {\Bbb D}_7+4 {\Bbb D}_{\infty} {\Bbb D}_7
               -2 {\Bbb D}_1 {\Bbb D}_6-2 {\Bbb D}_2 {\Bbb D}_5-2 {\Bbb D}_8 {\Bbb D}_{12}
               -2 {\Bbb D}_9 {\Bbb D}_{11},\\
  {\Bbb G}_8 &=-2 {\Bbb D}_4^2+6 {\Bbb D}_0 {\Bbb D}_8+4 {\Bbb D}_{\infty} {\Bbb D}_8
               -2 {\Bbb D}_3 {\Bbb D}_5-2 {\Bbb D}_6 {\Bbb D}_2-2 {\Bbb D}_{11} {\Bbb D}_{10}
               -2 {\Bbb D}_1 {\Bbb D}_7,\\
  {\Bbb G}_9 &=-{\Bbb D}_{11}^2+2 {\Bbb D}_0 {\Bbb D}_9+10 {\Bbb D}_{\infty} {\Bbb D}_9
               +2 {\Bbb D}_5 {\Bbb D}_4-2 {\Bbb D}_1 {\Bbb D}_8-4 {\Bbb D}_{10} {\Bbb D}_{12}
               -2 {\Bbb D}_3 {\Bbb D}_6,\\
  {\Bbb G}_{10} &=-{\Bbb D}_5^2+10 {\Bbb D}_0 {\Bbb D}_{10}-2 {\Bbb D}_{\infty} {\Bbb D}_{10}
               +2 {\Bbb D}_6 {\Bbb D}_4-2 {\Bbb D}_3 {\Bbb D}_7-4 {\Bbb D}_9 {\Bbb D}_1
               -2 {\Bbb D}_{12} {\Bbb D}_{11},\\
  {\Bbb G}_{11} &=-2 {\Bbb D}_{12}^2+6 {\Bbb D}_0 {\Bbb D}_{11}+4 {\Bbb D}_{\infty} {\Bbb D}_{11}
               -2 {\Bbb D}_9 {\Bbb D}_2-2 {\Bbb D}_5 {\Bbb D}_6-2 {\Bbb D}_7 {\Bbb D}_4
               -2 {\Bbb D}_3 {\Bbb D}_8,\\
  {\Bbb G}_{12} &=-{\Bbb D}_6^2+10 {\Bbb D}_0 {\Bbb D}_{12}-2 {\Bbb D}_{\infty} {\Bbb D}_{12}
               +2 {\Bbb D}_2 {\Bbb D}_{10}-2 {\Bbb D}_1 {\Bbb D}_{11}-4 {\Bbb D}_3 {\Bbb D}_9
               -2 {\Bbb D}_4 {\Bbb D}_8
\endaligned\right.\eqno{(1.10)}$$
are the senary sextic forms (sextic forms in six variables). Here,
$$\left\{\aligned
  {\Bbb D}_0 &=z_1 z_2 z_3,\\
  {\Bbb D}_1 &=2 z_2 z_3^2+z_2^2 z_6-z_4^2 z_5+z_1 z_5 z_6,\\
  {\Bbb D}_2 &=-z_6^3+z_2^2 z_4-2 z_2 z_5^2+z_1 z_4 z_5+3 z_3 z_5 z_6,\\
  {\Bbb D}_3 &=2 z_1 z_2^2+z_1^2 z_5-z_4 z_6^2+z_3 z_4 z_5,\\
  {\Bbb D}_4 &=-z_2^2 z_3+z_1 z_6^2-2 z_4^2 z_6-z_1 z_3 z_5,\\
  {\Bbb D}_5 &=-z_4^3+z_3^2 z_5-2 z_3 z_6^2+z_2 z_5 z_6+3 z_1 z_4 z_6,\\
  {\Bbb D}_6 &=-z_5^3+z_1^2 z_6-2 z_1 z_4^2+z_3 z_4 z_6+3 z_2 z_4 z_5,\\
  {\Bbb D}_7 &=-z_2^3+z_3 z_4^2-z_1 z_3 z_6-3 z_1 z_2 z_5+2 z_1^2 z_4,\\
  {\Bbb D}_8 &=-z_1^3+z_2 z_6^2-z_2 z_3 z_5-3 z_1 z_3 z_4+2 z_3^2 z_6,\\
  {\Bbb D}_9 &=2 z_1^2 z_3+z_3^2 z_4-z_5^2 z_6+z_2 z_4 z_6,\\
  {\Bbb D}_{10} &=-z_1 z_3^2+z_2 z_4^2-2 z_4 z_5^2-z_1 z_2 z_6,\\
  {\Bbb D}_{11} &=-z_3^3+z_1 z_5^2-z_1 z_2 z_4-3 z_2 z_3 z_6+2 z_2^2 z_5,\\
  {\Bbb D}_{12} &=-z_1^2 z_2+z_3 z_5^2-2 z_5 z_6^2-z_2 z_3 z_4,\\
  {\Bbb D}_{\infty}&=z_4 z_5 z_6
\endaligned\right.\eqno{(1.11)}$$
are the senary cubic forms (cubic forms in six variables).

{\bf Theorem 1.1. (Main Theorem 1)}. {\it The invariant
decomposition formula for the simple group $PSL(2, 13)$ of order
$1092$ is given as follows:
$$\left[\frac{\eta^2(z)}{\eta^2(13z)}\right]^5
 =\frac{\Phi_{12}(a_1(z), a_2(z), a_3(z), a_4(z), a_5(z), a_6(z))}
  {(a_1(z) a_2(z) a_3(z) a_4(z) a_5(z) a_6(z))^2},\eqno{(1.12)}$$
where $\Phi_{12}(z_1, z_2, z_3, z_4, z_5, z_6)$ is an invariant of
degree $12$ associated to $PSL(2, 13)$, and $(z_1 z_2 z_3 z_4 z_5
z_6)^2$ is invariant under the action of the image of a Borel
subgroup of $PSL(2, 13)$ $($i.e. a maximal subgroup of order $78$
of $PSL(2, 13)$$)$.}

  One of the significance of the invariant decomposition formulas (1.6),
(1.7) and (1.12) comes from that the Hauptmodul
$f(z)=(\eta(z)/\eta(pz))^{24/(p-1)}$ is a modular unit in the
rational function field of $X_0(p)$ for $p=5, 7, 13$. Its divisor
is given by $\text{div}(f)=\{ (0)-(\infty) \}$. Hence, (1.6),
(1.7) and (1.12) give some decomposition formulas of the modular
units which play an important role both in the function field
theory and applications to number theory (see \cite{KL} for more
details).

  The other significance of the invariant decomposition formulas (1.6),
(1.7) and (1.12) comes from the following fact: let $H$ be a
maximal subgroup of $PSL(2, p)$ for $p=5$, $7$ and $13$ with order
$10$, $21$ and $78$, respectively. Then
$$[\Gamma: \Gamma_0(p)]=[PSL(2, p): H]=p+1 \quad \text{for $p=5$, $7$
  and $13$, where $\Gamma=PSL(2, {\Bbb Z})$}.$$

  Perhaps the most significance of the invariant decomposition formulas
(1.6), (1.7) and (1.12) is the relation with the Monster simple
group ${\Bbb M}$ and monstrous moonshine. In fact, elements of
type $5B$, $7B$, $13B$ in the Monster correspond to elements in
the conjugacy classes $5A$, $7A$, $13A$ of the Conway group $Co_1$
with cycle shapes $5^6 1^{-6}$, $7^4 1^{-4}$, $13^2 1^{-2}$. The
modular groups corresponding to these elements of the Monster are
$\Gamma_0(5)$, $\Gamma_0(7)$ and $\Gamma_0(13)$ (see \cite{Bor},
\cite{CN} and \cite{CC}). Moreover, the Hauptmoduln
$(\eta(z)/\eta(5z))^6$, $(\eta(z)/\eta(7z))^4$ and
$(\eta(z)/\eta(13z))^2$ correspond exactly to the cycle shapes
$5^6 1^{-6}$, $7^4 1^{-4}$, $13^2 1^{-2}$, respectively (see
Appendix for more details).

  In \cite{K1}, Klein obtained the modular equation of degree $14$
for $\Gamma_0(13)$ (see \cite{CY}, p.267 for the modular relation
for $\Gamma_0(13)$). Combining with Theorem 1.1, we give a new
expression of the classical $j$-function in terms of theta
constants associated with $\Gamma(13)$.

{\bf Corollary 1.2}. {\it The following formulas hold:
$$j(z)=\frac{(\tau^2+5\tau+13)(\tau^4+247 \tau^3+3380 \tau^2+15379 \tau+28561)^3}{\tau^{13}}\eqno{(1.13)}$$
and
$$j(13z)=\frac{(\tau^2+5 \tau+13)(\tau^4+7 \tau^3+20 \tau^2+19 \tau+1)^3}{\tau},\eqno{(1.14)}$$
where the Hauptmodul
$$\tau=\left(\frac{\eta(z)}{\eta(13z)}\right)^2
   =\root 5 \of{\frac{\Phi_{12}(a_1(z), a_2(z), a_3(z), a_4(z), a_5(z), a_6(z))}
  {(a_1(z) a_2(z) a_3(z) a_4(z) a_5(z) a_6(z))^2}}.\eqno{(1.15)}$$}

  Note that in the right hand-side of modular equations (1.13) and (1.14),
we have the following factorizations over ${\Bbb Q}(\sqrt{13})$:
$$\aligned
 &\tau^4+247 \tau^3+3380 \tau^2+15379 \tau+28561\\
=&\left(\tau^2+\frac{247+65 \sqrt{13}}{2} \tau+\frac{1859+507
  \sqrt{13}}{2}\right)\left(\tau^2+\frac{247-65 \sqrt{13}}{2}
  \tau+\frac{1859-507 \sqrt{13}}{2}\right)
\endaligned$$
and
$$\aligned
 &\tau^4+7 \tau^3+20 \tau^2+19 \tau+1\\
=&\left(\tau^2+\frac{7+\sqrt{13}}{2} \tau+\frac{11+3
  \sqrt{13}}{2}\right)\left(\tau^2+\frac{7-\sqrt{13}}{2} \tau+\frac{11-3 \sqrt{13}}{2}\right).
\endaligned$$
In \cite{De1}, Deligne gave a modular interpretation of why such
factorizations exist. More generally, in \cite{De2}, Deligne
showed that for $p=3, 5, 7, 13$, the corresponding modular
equations of degree $4$, $6$, $8$, $14$ have two relatively prime
conjugate factors over some real quadratic fields. All of these
factorizations have nice geometric interpretations. This is the
other significance of the invariant decomposition formulas (1.6),
(1.7) and (1.12) since they give the expressions of the Hauptmodul
$\tau=(\eta(z)/\eta(pz))^{24/(p-1)}$ for $p=5, 7, 13$.

  In his notebooks (see \cite{R5}, p.326 and \cite{R6}, p.244), Ramanujan
(see also \cite{Ev} or \cite{B1}, p.373-375) obtained the
following modular equation of degree $13$:
$$1+G_1(z) G_5(z)+G_2(z) G_3(z)+G_4(z) G_6(z)=-\frac{\eta^2(z/13)}{\eta^2(z)}\eqno{(1.16)}$$
with
$$G_1(z)G_5(z)G_2(z)G_3(z)G_4(z)G_6(z)=-1,\eqno{(1.17)}$$
where
$$G_m(z):=G_{m, p}(z):=(-1)^m q^{m(3m-p)/(2p^2)} \frac{f(-q^{2m/p}, -q^{1-2m/p})}{f(-q^{m/p}, -q^{1-m/p})}.$$
Here, $m$ is a positive integer, $p=13$, $q=e^{2 \pi i z}$, and
Ramanujan's general theta function $f(a, b)$ is given by
$$f(a, b):=\sum_{n=-\infty}^{\infty} a^{n(n+1)/2} b^{n(n-1)/2}, \quad |ab|<1,$$
We call (1.16) the standard modular equation of degree $13$. In
contrast with it, we find the following invariant decomposition
formula (exotic modular equation) which has the same form as
(1.16), but with different kinds of modular parametrizations. Let
us define the following senary quadratic forms (quadratic forms in
six variables):
$$\left\{\aligned
  {\Bbb A}_0 &=z_1 z_4+z_2 z_5+z_3 z_6,\\
  {\Bbb A}_1 &=z_1^2-2 z_3 z_4,\\
  {\Bbb A}_2 &=-z_5^2-2 z_2 z_4,\\
  {\Bbb A}_3 &=z_2^2-2 z_1 z_5,\\
  {\Bbb A}_4 &=z_3^2-2 z_2 z_6,\\
  {\Bbb A}_5 &=-z_4^2-2 z_1 z_6,\\
  {\Bbb A}_6 &=-z_6^2-2 z_3 z_5.
\endaligned\right.\eqno{(1.18)}$$

{\bf Theorem 1.3. (Main Theorem 2)}. {\it The following invariant
decomposition formula $($exotic modular equation$)$ holds:
$$\frac{\Psi_2(a_1(z), a_2(z), a_3(z), a_4(z), a_5(z), a_6(z))}
  {{\Bbb A}_0(a_1(z), a_2(z), a_3(z), a_4(z), a_5(z), a_6(z))^2}=0,\eqno{(1.19)}$$
where the quadric
$$\Psi_2(z_1, z_2, z_3, z_4, z_5, z_6):=
 {\Bbb A}_0^2+{\Bbb A}_1 {\Bbb A}_5+{\Bbb A}_2 {\Bbb A}_3+{\Bbb A}_4 {\Bbb A}_6\eqno{(1.20)}$$
is an invariant associated to $PSL(2, 13)$ and ${\Bbb A}_0^2$ is
invariant under the action of the image of a Borel subgroup of
$PSL(2, 13)$ $($i.e. a maximal subgroup of order $78$ of $PSL(2,
13)$$)$.}

  Note that (1.19) has the same form as (1.16). However, in contrast with (1.17),
it can be proved that (see section 3, (3.22))
$$\prod_{j=1}^{6} \frac{{\Bbb A}_j(a_1(z), a_2(z), a_3(z), a_4(z), a_5(z), a_6(z))}
       {{\Bbb A}_0(a_1(z), a_2(z), a_3(z), a_4(z), a_5(z), a_6(z))} \neq -1.\eqno{(1.21)}$$
Moreover, the invariant quadric (1.20) is closely related to the
exceptional Lie group $G_2$ (see \cite{Y} for more details).

  In \cite{K}, Klein investigated the modular equation of degree
six in connection with the Jacobian equation of degree six, which
can be used to solve the general quintic equation. In \cite{K2}
and \cite{K3}, he studied the modular equation of degree eight in
connection with the Jacobian equation of degree eight, which can
be used to solve the algebraic equation of degree eight with
Galois group $PSL(2, 7)$. However, he did not investigate the last
case, the modular equation of degree fourteen. In this paper, we
find the Jacobian equation of degree fourteen, which corresponds
to the modular equation of degree fourteen. Surprisingly, we find
the other exotic modular equation of degree fourteen, which does
not correspond to the Jacobian equation. It plays an central role
in the proof of Theorem 1.1. As an application, we obtain the
following quartic four-fold $\Phi_4(z_1, z_2, z_3, z_4, z_5,
z_6)=0$, where
$$\Phi_4:=(z_3 z_4^3+z_1 z_5^3+z_2 z_6^3)-(z_6 z_1^3+z_4 z_2^3+z_5 z_3^3)+
  3(z_1 z_2 z_4 z_5+z_2 z_3 z_5 z_6+z_3 z_1 z_6 z_4),\eqno{(1.22)}$$
which is just the quadric (1.20) up to a constant. It is a
higher-dimensional counterpart of the Klein quartic curve (see
\cite{K2}) and the Klein cubic threefold (see \cite{K4}). Its
significance comes from that it gives the geometry of the
classical modular curve $X(13)$ of genus $50$:

{\bf Corollary 1.4}. {\it The coordinates $(a_1(z), a_2(z),
a_3(z), a_4(z), a_5(z), a_6(z))$ map $X(13)$ into the quartic
four-fold $\Phi_4(z_1, z_2, z_3, z_4, z_5, z_6)=0$ in ${\Bbb C}
{\Bbb P}^5$.}

  This paper consists of four sections. In section two, we give
a six-dimensional representation of $PSL(2, 13)$ defined over
${\Bbb Q}(e^{\frac{2 \pi i}{13}})$ and the transformation formulas
for theta constants associated with $\Gamma(13)$. In section
three, we give a seven-dimensional representation of $PSL(2, 13)$
which is deduced from our six-dimensional representation. As an
application, we obtain the exotic modular equation of degree $13$
and give the geometry of the modular curve $X(13)$. Thus, we give
the proof of Theorem 1.3. In section four, we give a
fourteen-dimensional representation of $PSL(2, 13)$ which is also
deduced from our six-dimensional representation. By this
representation, we find the exotic modular equation of degree
fourteen. Using it, we give the proof of Theorem 1.1.

\vskip 0.5 cm

\flushpar
{\bf Acknowedgements}: I thank Pierre Deligne for his
comments on an earlier version of this paper.

\vskip 0.5 cm

\centerline{\bf 2. Six-dimensional representations of $PSL(2, 13)$
                   and}

\centerline{\bf transformation formulas for theta constants}

\vskip 0.5 cm

  Let us recall the Weil representation of $SL_2({\Bbb F}_p)$
(see \cite{W}) which is modelled on the space $L^2({\Bbb F}_p)$ of
all complex valued functions on ${\Bbb F}_p$ where $p$ is a prime.
Denote by $L^2({\Bbb F}_p)$ the $p$-dimensional complex vector
space of all (square-integrable) complex valued functions on
${\Bbb F}_p$ with respect to counting measure, that is, all
functions from ${\Bbb F}_p$ to the complex numbers. We can
decompose $L^2({\Bbb F}_p)$ as the direct sum of the space $V^{+}$
of even functions and the space $V^{-}$ of odd functions. The
space $V^{-}$ has dimension $(p-1)/2$ and its associated
projective space ${\Bbb P}(V^{-})$ has dimension $(p-3)/2$. In
particular, for $p=13$, the Weil representation for $SL(2, 13)$ is
given as follows:
$$s=\left(\matrix 0 & -1\\ 1 & 0 \endmatrix\right), \quad
  t=\left(\matrix 1 & 1\\ 0 & 1 \endmatrix\right), \quad
  h=\left(\matrix 7 & 0\\ 0 & 2 \endmatrix\right).$$

  In this section, we will study the six-dimensional representation
of the simple group $PSL(2, 13)$ of order $1092$, which acts on
the five-dimensional projective space
$${\Bbb P}^5=\{ (z_1, z_2, z_3, z_4, z_5, z_6): z_i \in {\Bbb C} \quad (i=1, 2, 3, 4, 5, 6) \}.$$

  Let $\zeta=\exp(2 \pi i/13)$ and
$$\left\{\aligned
  &\theta_1=\zeta+\zeta^3+\zeta^9,\\
  &\theta_2=\zeta^2+\zeta^6+\zeta^5,\\
  &\theta_3=\zeta^4+\zeta^{12}+\zeta^{10},\\
  &\theta_4=\zeta^8+\zeta^{11}+\zeta^7.
\endaligned\right.\eqno{(2.1)}$$
We find that
$$\left\{\aligned
  &\theta_1+\theta_2+\theta_3+\theta_4=-1,\\
  &\theta_1 \theta_2+\theta_1 \theta_3+\theta_1 \theta_4+\theta_2 \theta_3+\theta_2 \theta_4+\theta_3 \theta_4=2,\\
  &\theta_1 \theta_2 \theta_3+\theta_1 \theta_2 \theta_4+\theta_1 \theta_3 \theta_4+\theta_2 \theta_3 \theta_4=4,\\
  &\theta_1 \theta_2 \theta_3 \theta_4=3.
\endaligned\right.$$
Hence, $\theta_1$, $\theta_2$, $\theta_3$ and $\theta_4$ satisfy
the quartic equation
$$z^4+z^3+2 z^2-4z+3=0,$$
which can be decomposed as two quadratic equations
$$\left(z^2+\frac{1+\sqrt{13}}{2} z+\frac{5+\sqrt{13}}{2}\right)
  \left(z^2+\frac{1-\sqrt{13}}{2} z+\frac{5-\sqrt{13}}{2}\right)=0$$
over the real quadratic field ${\Bbb Q}(\sqrt{13})$. Therefore,
the four roots are given as follows:
$$z_{1, 2}=\frac{1}{2}\left(-\frac{1+\sqrt{13}}{2} \pm \sqrt{\frac{-13-3 \sqrt{13}}{2}}\right)
          =\frac{1}{4} \left(-1-\sqrt{13} \pm \sqrt{-26-6 \sqrt{13}}\right),$$
$$z_{3, 4}=\frac{1}{2}\left(-\frac{1-\sqrt{13}}{2} \pm \sqrt{\frac{-13+3 \sqrt{13}}{2}}\right)
          =\frac{1}{4} \left(-1+\sqrt{13} \pm \sqrt{-26+6
\sqrt{13}}\right).$$ Note that
$$\text{Re}(\theta_1)=\cos \frac{2 \pi}{13}+\cos \frac{6 \pi}{13}-\cos \frac{5 \pi}{13}>0, \quad
  \text{Im}(\theta_1)=\sin \frac{2 \pi}{13}+\sin \frac{6 \pi}{13}-\sin \frac{5 \pi}{13}>0.$$
We have
$$\theta_1=\frac{1}{4} \left(-1+\sqrt{13}+\sqrt{-26+6 \sqrt{13}}\right).$$
Similarly,
$$\text{Re}(\theta_2)<0, \quad \text{Im}(\theta_2)>0, \quad
  \theta_2=\frac{1}{4} \left(-1-\sqrt{13}+\sqrt{-26-6 \sqrt{13}}\right).$$
$$\text{Re}(\theta_3)>0, \quad \text{Im}(\theta_3)<0, \quad
  \theta_3=\frac{1}{4} \left(-1+\sqrt{13}-\sqrt{-26+6 \sqrt{13}}\right).$$
$$\text{Re}(\theta_4)<0, \quad \text{Im}(\theta_4)<0, \quad
  \theta_4=\frac{1}{4} \left(-1-\sqrt{13}-\sqrt{-26-6 \sqrt{13}}\right).$$
Moreover, we find that
$$\left\{\aligned
  \theta_1+\theta_3+\theta_2+\theta_4 &=-1,\\
  \theta_1+\theta_3-\theta_2-\theta_4 &=\sqrt{13},\\
  \theta_1-\theta_3-\theta_2+\theta_4 &=-\sqrt{-13+2 \sqrt{13}},\\
  \theta_1-\theta_3+\theta_2-\theta_4 &=\sqrt{-13-2 \sqrt{13}}.
\endaligned\right.\eqno{(2.2)}$$

  Let
$$M=\left(\matrix
    \zeta-\zeta^{12} & \zeta^3-\zeta^{10} & \zeta^9-\zeta^4\\
    \zeta^3-\zeta^{10} & \zeta^9-\zeta^4 & \zeta-\zeta^{12}\\
    \zeta^9-\zeta^4 & \zeta-\zeta^{12} & \zeta^3-\zeta^{10}
   \endmatrix\right), \quad
  N=\left(\matrix
    \zeta^5-\zeta^8 & \zeta^2-\zeta^{11} & \zeta^6-\zeta^7\\
    \zeta^2-\zeta^{11} & \zeta^6-\zeta^7 & \zeta^5-\zeta^8\\
    \zeta^6-\zeta^7 & \zeta^5-\zeta^8 & \zeta^2-\zeta^{11}
   \endmatrix\right).\eqno{(2.3)}$$
Then $MN=NM=-\sqrt{13} I$ and $M^2+N^2=-13 I$. Put
$$S=-\frac{1}{\sqrt{13}} \left(\matrix
    -M & N\\
     N & M
    \endmatrix\right).\eqno{(2.4)}$$
Then $S^2=I$. In fact,
$$S=-\frac{2 i}{\sqrt{13}} \left(\matrix
    -\sin \frac{2 \pi}{13} & -\sin \frac{6 \pi}{13} &  \sin \frac{5 \pi}{13} &
     \sin \frac{3 \pi}{13} &  \sin \frac{4 \pi}{13} &  \sin \frac{\pi}{13}\\
    -\sin \frac{6 \pi}{13} &  \sin \frac{5 \pi}{13} & -\sin \frac{2 \pi}{13} &
     \sin \frac{4 \pi}{13} &  \sin \frac{\pi}{13}   &  \sin \frac{3 \pi}{13}\\
     \sin \frac{5 \pi}{13} & -\sin \frac{2 \pi}{13} & -\sin \frac{6 \pi}{13} &
     \sin \frac{\pi}{13}   &  \sin \frac{3 \pi}{13} &  \sin \frac{4 \pi}{13}\\
     \sin \frac{3 \pi}{13} &  \sin \frac{4 \pi}{13} &  \sin \frac{\pi}{13} &
     \sin \frac{2 \pi}{13} &  \sin \frac{6 \pi}{13} & -\sin \frac{5 \pi}{13}\\
     \sin \frac{4 \pi}{13} &  \sin \frac{\pi}{13}   &  \sin \frac{3 \pi}{13} &
     \sin \frac{6 \pi}{13} & -\sin \frac{5 \pi}{13} &  \sin \frac{2 \pi}{13}\\
     \sin \frac{\pi}{13}   &  \sin \frac{3 \pi}{13} &  \sin \frac{4 \pi}{13} &
    -\sin \frac{5 \pi}{13} &  \sin \frac{2 \pi}{13} &  \sin \frac{6 \pi}{13}
   \endmatrix\right).\eqno{(2.5)}$$
Note that
$$\frac{\sin \frac{2 \pi}{13} \sin \frac{5 \pi}{13} \sin \frac{6 \pi}{13}}
       {\sin \frac{\pi}{13} \sin \frac{3 \pi}{13} \sin \frac{4 \pi}{13}}
 =\frac{-\sqrt{\frac{-13-3 \sqrt{13}}{2}}}{-\sqrt{\frac{-13+3 \sqrt{13}}{2}}}
 =\frac{3+\sqrt{13}}{2},$$
which is a fundamental unit of ${\Bbb Q}(\sqrt{13})$. Let
$$T=\text{diag}(\zeta^7, \zeta^{11}, \zeta^8, \zeta^6, \zeta^2, \zeta^5).\eqno{(2.6)}$$

{\bf Theorem 2.1.} {\it Let $G=\langle S, T \rangle$. Then $G
\cong PSL(2, 13)$.}

{\it Proof}. We have
$$ST=-\frac{1}{\sqrt{13}} \left(\matrix
  \zeta^6-\zeta^8 & \zeta^8-\zeta & \zeta^{12}-\zeta^4 & \zeta^{11}-\zeta & \zeta^4-1 & \zeta^{11}-\zeta^{12}\\
  \zeta^4-\zeta^{10} & \zeta^2-\zeta^7 & \zeta^7-\zeta^9 & \zeta^8-\zeta^4 & \zeta^8-\zeta^9 & \zeta^{10}-1\\
  \zeta^{11}-\zeta^3 & \zeta^{10}-\zeta^{12} & \zeta^5-\zeta^{11} & \zeta^{12}-1 & \zeta^7-\zeta^{10} & \zeta^7-\zeta^3\\
  \zeta^{12}-\zeta^2 & 1-\zeta^9 & \zeta-\zeta^2 & \zeta^7-\zeta^5 & \zeta^5-\zeta^{12} & \zeta-\zeta^9\\
  \zeta^9-\zeta^5 & \zeta^4-\zeta^5 & 1-\zeta^3 & \zeta^9-\zeta^3 & \zeta^{11}-\zeta^6 & \zeta^6-\zeta^4\\
  1-\zeta & \zeta^3-\zeta^6 & \zeta^{10}-\zeta^6 & \zeta^2-\zeta^{10} & \zeta^3-\zeta & \zeta^8-\zeta^2
\endmatrix\right).\eqno{(2.7)}$$
On the other hand,
$$\aligned
 &(ST)^{-1}=T^{-1} S\\
=&-\frac{1}{\sqrt{13}} \left(\matrix
  \zeta^5-\zeta^7 & \zeta^3-\zeta^9 & \zeta^{10}-\zeta^2 & \zeta^{11}-\zeta & \zeta^8-\zeta^4 & \zeta^{12}-1\\
  \zeta^{12}-\zeta^5 & \zeta^6-\zeta^{11} & \zeta-\zeta^3 & \zeta^4-1 & \zeta^8-\zeta^9 & \zeta^7-\zeta^{10}\\
  \zeta^9-\zeta & \zeta^4-\zeta^6 & \zeta^2-\zeta^8 & \zeta^{11}-\zeta^{12} & \zeta^{10}-1 & \zeta^7-\zeta^3\\
  \zeta^{12}-\zeta^2 & \zeta^9-\zeta^5 & 1-\zeta & \zeta^8-\zeta^6 & \zeta^{10}-\zeta^4 & \zeta^3-\zeta^{11}\\
  1-\zeta^9 & \zeta^4-\zeta^5 & \zeta^3-\zeta^6 & \zeta-\zeta^8 & \zeta^7-\zeta^2 & \zeta^{12}-\zeta^{10}\\
  \zeta-\zeta^2 & 1-\zeta^3 & \zeta^{10}-\zeta^6 & \zeta^4-\zeta^{12} & \zeta^9-\zeta^7 & \zeta^{11}-\zeta^5
\endmatrix\right).
\endaligned$$
We will calculate $(ST)^2=\frac{1}{13} (a_{ij})$. Without loss of
generality, we will compute $a_{1i}$, $i=1, 2, 3, 4, 5, 6$. Here,
$$a_{11}=-2 \zeta-2 \zeta^2+2 \zeta^3-2 \zeta^9+2 \zeta^{10}+2\zeta^{11}-\zeta^5+\zeta^7.$$
By the identity
$$\sqrt{13}=\zeta+\zeta^{12}+\zeta^3+\zeta^{10}+\zeta^9+\zeta^4-\zeta^5-\zeta^8-\zeta^2-\zeta^{11}-\zeta^6-\zeta^7,\eqno{(2.8)}$$
we have
$$(\zeta^5-\zeta^7) \sqrt{13}
 =2 \zeta+2 \zeta^2-2 \zeta^3+2 \zeta^9-2 \zeta^{10}-2 \zeta^{11}+\zeta^5-\zeta^7.$$
Hence, $a_{11}=-(\zeta^5-\zeta^7) \sqrt{13}$. Similarly,
$$\aligned
  a_{12}&=-2 \zeta^2-2 \zeta^4+2 \zeta^5-2 \zeta^7+2 \zeta^8+2 \zeta^{10}-\zeta^3+\zeta^9
 =-(\zeta^3-\zeta^9) \sqrt{13},\\
  a_{13}&=-2+2 \zeta^3+2 \zeta^5-2 \zeta^7-2 \zeta^9+2 \zeta^{12}+\zeta^2-\zeta^{10}
 =-(\zeta^{10}-\zeta^2) \sqrt{13},\\
  a_{14}&=2+2 \zeta^4+2 \zeta^5-2 \zeta^7-2 \zeta^8-2 \zeta^{12}+\zeta^{11}-\zeta
 =-(\zeta^{11}-\zeta) \sqrt{13},\\
  a_{15}&=2+2 \zeta+2 \zeta^3-2 \zeta^9-2 \zeta^{11}-2 \zeta^{12}+\zeta^8-\zeta^4
 =-(\zeta^8-\zeta^4) \sqrt{13},\\
  a_{16}&=2 \zeta-2 \zeta^2+2 \zeta^4-2 \zeta^8+2 \zeta^{10}-2 \zeta^{11}+\zeta^{12}-1
 =-(\zeta^{12}-1) \sqrt{13}.
\endaligned$$
The other terms can be calculated in the same way. In conclusion,
we find that $(ST)^2=(ST)^{-1}$, i.e., $(ST)^3=1$. Hence, we have
$$S^2=T^{13}=(ST)^3=1.\eqno{(2.9)}$$
Let $u=ST$ and $v=S$. Then $uv=STS$. Hence,
$$u^3=v^2=(uv)^{13}=1.\eqno{(2.10)}$$
According to \cite{S}, put
$$(2, 3, n; p):=\langle u, v: u^3=v^2=(uv)^n=(u^{-1} v^{-1} uv)^p=1 \rangle.$$
Let $P=(uv)^{-1}$ and $Q=(uv)^2 u$. Then $u=P^2 Q$ and $v=P^3 Q$.
In \cite{S}, Sinkov proved the following:

{\bf Theorem 2.2.} (see \cite{S}). {\it Two operators of periods
$2$ and $3$ generate the simple group of order $1092$ if and only
if they satisfy one of the following sets of independent
relations:
$$\aligned
  A: \quad &(2, 3, 7; 6),\\
  B: \quad &(2, 3, 7; 7),\\
  C: \quad &(2, 3, 7); (Q^2 P^6)^3=1,\\
  D: \quad &(2, 3, 13); (Q^3 P^4)^3=1.
\endaligned$$}

  In our case, $P=S T^{-1} S$ and $Q=S T^3$.
We have
$$Q^3=-\frac{1}{\sqrt{13}} \left(\matrix
  \zeta^{11}-\zeta & \zeta^{12}-\zeta^8 & 1-\zeta & \zeta^6-\zeta^4 & \zeta^4-\zeta^{11} & 1-\zeta^8\\
  1-\zeta^9 & \zeta^8-\zeta^9 & \zeta^4-\zeta^7 & 1-\zeta^7 & \zeta^2-\zeta^{10} & \zeta^{10}-\zeta^8\\
  \zeta^{10}-\zeta^{11} & 1-\zeta^3 & \zeta^7-\zeta^3 & \zeta^{12}-\zeta^7 & 1-\zeta^{11} & \zeta^5-\zeta^{12}\\
  \zeta^9-\zeta^7 & \zeta^2-\zeta^9 & \zeta^5-1 & \zeta^2-\zeta^{12} & \zeta-\zeta^5 & 1-\zeta^{12}\\
  \zeta^6-1 & \zeta^3-\zeta^{11} & \zeta^5-\zeta^3 & 1-\zeta^4 & \zeta^5-\zeta^4 & \zeta^9-\zeta^6\\
  \zeta^6-\zeta & \zeta^2-1 & \zeta-\zeta^8 & \zeta^3-\zeta^2 & 1-\zeta^{10} & \zeta^6-\zeta^{10}
\endmatrix\right).$$
$$P^4=-\frac{1}{\sqrt{13}} \left(\matrix
  \zeta^7-1 & \zeta^2-\zeta^7 & \zeta^6-\zeta^8 & \zeta^2-\zeta^{11} & \zeta^5-\zeta^6 & \zeta^8-\zeta^{11}\\
  \zeta^2-\zeta^7 & \zeta^{11}-1 & \zeta^5-\zeta^{11} & \zeta^7-\zeta^8 & \zeta^5-\zeta^8 & \zeta^6-\zeta^2\\
  \zeta^6-\zeta^8 & \zeta^5-\zeta^{11} & \zeta^8-1 & \zeta^2-\zeta^5 & \zeta^{11}-\zeta^7 & \zeta^6-\zeta^7\\
  \zeta^2-\zeta^{11} & \zeta^7-\zeta^8 & \zeta^2-\zeta^5 & \zeta^6-1 & \zeta^{11}-\zeta^6 & \zeta^7-\zeta^5\\
  \zeta^5-\zeta^6 & \zeta^5-\zeta^8 & \zeta^{11}-\zeta^7 & \zeta^{11}-\zeta^6 & \zeta^2-1 & \zeta^8-\zeta^2\\
  \zeta^8-\zeta^{11} & \zeta^6-\zeta^2 & \zeta^6-\zeta^7 & \zeta^7-\zeta^5 & \zeta^8-\zeta^2 & \zeta^5-1
\endmatrix\right).$$
$$Q^3 P^4=-\frac{1}{\sqrt{13}} \left(\matrix
  \zeta^7-\zeta^5 & \zeta^2-\zeta^9 & \zeta^{10}-\zeta^5 & \zeta^6-\zeta^3 & \zeta^3-\zeta^7 & \zeta^{10}-\zeta^9\\
  \zeta^{12}-\zeta^6 & \zeta^{11}-\zeta^6 & \zeta^5-\zeta^3 & \zeta^{12}-\zeta^3 & \zeta^2-\zeta & \zeta-\zeta^{11}\\
  \zeta^6-\zeta & \zeta^4-\zeta^2 & \zeta^8-\zeta^2 & \zeta^9-\zeta^8 & \zeta^4-\zeta & \zeta^5-\zeta^9\\
  \zeta^{10}-\zeta^7 & \zeta^6-\zeta^{10} & \zeta^4-\zeta^3 & \zeta^6-\zeta^8 & \zeta^{11}-\zeta^4 & \zeta^3-\zeta^8\\
  \zeta^{10}-\zeta & \zeta^{12}-\zeta^{11} & \zeta^2-\zeta^{12} & \zeta-\zeta^7 & \zeta^2-\zeta^7 & \zeta^8-\zeta^{10}\\
  \zeta^5-\zeta^4 & \zeta^{12}-\zeta^9 & \zeta^4-\zeta^8 & \zeta^7-\zeta^{12} & \zeta^9-\zeta^{11} & \zeta^5-\zeta^{11}
\endmatrix\right),$$
$$(Q^3 P^4)^2=-\frac{1}{\sqrt{13}} \left(\matrix
  \zeta^8-\zeta^6 &  \zeta^7-\zeta & \zeta^{12}-\zeta^7 & \zeta^6-\zeta^3 & \zeta^{12}-\zeta^3 & \zeta^9-\zeta^8\\
  \zeta^4-\zeta^{11} & \zeta^7-\zeta^2 & \zeta^{11}-\zeta^9 & \zeta^3-\zeta^7 & \zeta^2-\zeta & \zeta^4-\zeta\\
  \zeta^8-\zeta^3 & \zeta^{10}-\zeta^8 & \zeta^{11}-\zeta^5 & \zeta^{10}-\zeta^9 & \zeta-\zeta^{11} & \zeta^5-\zeta^9\\
  \zeta^{10}-\zeta^7 & \zeta^{10}-\zeta & \zeta^5-\zeta^4 & \zeta^5-\zeta^7 & \zeta^6-\zeta^{12} & \zeta-\zeta^6\\
  \zeta^6-\zeta^{10} & \zeta^{12}-\zeta^{11} & \zeta^{12}-\zeta^9 & \zeta^9-\zeta^2 & \zeta^6-\zeta^{11} & \zeta^2-\zeta^4\\
  \zeta^4-\zeta^3 & \zeta^2-\zeta^{12} & \zeta^4-\zeta^8 & \zeta^5-\zeta^{10} & \zeta^3-\zeta^5 & \zeta^2-\zeta^8
\endmatrix\right),$$
and $(Q^3 P^4)^3=-I$. Note that in the projective coordinates,
this means that $(Q^3 P^4)^3=1$. Hence, we prove that the elements
$u$ and $v$ above satisfy the following relations: $(2, 3, 13)$
and $(Q^3 P^4)^3=1$, which is a presentation for the simple group
$PSL(2, 13)$ of order $1092$ by Theorem 3.2. Since the group is
simple and the generating matrices are non-trivial, we must have
$G=\langle u, v \rangle \cong PSL(2, 13)$. Hence, $\langle P, Q
\rangle=\langle S, T \rangle=G$. This completes the proof of
Theorem 2.1.

\flushpar $\qquad \qquad \qquad \qquad \qquad \qquad \qquad \qquad
\qquad \qquad \qquad \qquad \qquad \qquad \qquad \qquad \qquad
\qquad \quad \boxed{}$

  Recall that the theta functions with characteristic
$\left[\matrix \epsilon\\ \epsilon^{\prime} \endmatrix\right] \in
{\Bbb R}^2$ is defined by the following series which converges
uniformly and absolutely on compact subsets of ${\Bbb C} \times
{\Bbb H}$ (see \cite{FK}, p. 73):
$$\theta \left[\matrix \epsilon\\ \epsilon^{\prime} \endmatrix\right] (z, \tau)
 =\sum_{n \in {\Bbb Z}} \exp \left\{2 \pi i \left[\frac{1}{2}
  \left(n+\frac{\epsilon}{2}\right)^2 \tau+\left(n+\frac{\epsilon}{2}\right)
  \left(z+\frac{\epsilon^{\prime}}{2}\right)\right]\right\}.$$
The modified theta constants is defined by (see \cite{FK}, p. 215)
$$\varphi_l(\tau)=\theta [\chi_l](0, k \tau),$$
where the characteristic $\chi_l=\left[\matrix \frac{2l+1}{k}\\ 1
\endmatrix\right]$, $l=0, \cdots, \frac{k-3}{2}$, for odd $k$ and
$\chi_l=\left[\matrix \frac{2l}{k}\\ 0 \endmatrix\right]$, $l=0,
\cdots, \frac{k}{2}$, for even $k$. We have the following:

{\bf Theorem 2.3.} (see \cite{FK}, p. 236). {\it For each odd
integer $k \geq 5$, the map
$$\Phi: \tau \mapsto (\varphi_0(\tau), \varphi_1(\tau), \cdots,
  \varphi_{\frac{k-5}{2}}(\tau), \varphi_{\frac{k-3}{2}}(\tau))$$
from ${\Bbb H} \cup {\Bbb Q} \cup \{ \infty \}$ to ${\Bbb
C}^{\frac{k-1}{2}}$, defines a holomorphic mapping from
$\overline{{\Bbb H}/\Gamma(k)}$ into ${\Bbb C} {\Bbb
P}^{\frac{k-3}{2}}$.}

  In our case, $k=13$, the map
$$\Phi: \tau \mapsto (\varphi_0(\tau), \varphi_1(\tau), \varphi_2(\tau),
  \varphi_3(\tau), \varphi_4(\tau), \varphi_5(\tau))$$
gives a holomorphic mapping from the modular curve
$X(13)=\overline{{\Bbb H}/\Gamma(13)}$ into ${\Bbb C} {\Bbb P}^5$,
which corresponds to our six-dimensional representation, i.e., up
to the constants, $z_1, \cdots, z_6$ are just modular forms
$\varphi_0(\tau), \cdots, \varphi_5(\tau)$.

  Let
$$\left\{\aligned
  a_1(z) &:=e^{-\frac{11 \pi i}{26}} \theta \left[\matrix \frac{11}{13}\\ 1 \endmatrix\right](0, 13z)
           =q^{\frac{121}{104}} \sum_{n \in {\Bbb Z}} (-1)^n q^{\frac{1}{2}(13n^2+11n)},\\
  a_2(z) &:=e^{-\frac{7 \pi i}{26}} \theta \left[\matrix \frac{7}{13}\\ 1 \endmatrix\right](0, 13z)
           =q^{\frac{49}{104}} \sum_{n \in {\Bbb Z}} (-1)^n q^{\frac{1}{2}(13n^2+7n)},\\
  a_3(z) &:=e^{-\frac{5 \pi i}{26}} \theta \left[\matrix \frac{5}{13}\\ 1 \endmatrix\right](0, 13z)
           =q^{\frac{25}{104}} \sum_{n \in {\Bbb Z}} (-1)^n q^{\frac{1}{2}(13n^2+5n)},\\
  a_4(z) &:=-e^{-\frac{3 \pi i}{26}} \theta \left[\matrix \frac{3}{13}\\ 1 \endmatrix\right](0, 13z)
           =-q^{\frac{9}{104}} \sum_{n \in {\Bbb Z}} (-1)^n q^{\frac{1}{2}(13n^2+3n)},\\
  a_5(z) &:=e^{-\frac{9 \pi i}{26}} \theta \left[\matrix \frac{9}{13}\\ 1 \endmatrix\right](0, 13z)
           =q^{\frac{81}{104}} \sum_{n \in {\Bbb Z}} (-1)^n q^{\frac{1}{2}(13n^2+9n)},\\
  a_6(z) &:=e^{-\frac{\pi i}{26}} \theta \left[\matrix \frac{1}{13}\\ 1 \endmatrix\right](0, 13z)
           =q^{\frac{1}{104}} \sum_{n \in {\Bbb Z}} (-1)^n q^{\frac{1}{2}(13n^2+n)}.
\endaligned\right.\eqno{(2.11)}$$
We use the standard notation $(a)_{\infty}:=(a;
q):=\prod_{k=0}^{\infty} (1-aq^k)$ for $a \in {\Bbb C}^{\times}$
and $|q|<1$ so that $\eta(z)=q^{1/24}(q; q)$. By the Jacobi triple
product identity, we have
$$\theta \left[\matrix \frac{l}{k}\\ 1 \endmatrix\right](0, kz)
 =\exp \left(\frac{\pi i l}{2k}\right) q^{\frac{l^2}{8k}}
  (q^{\frac{k-l}{2}}; q^k)(q^{\frac{k+l}{2}}; q^k)(q^k; q^k),$$
where $k$ is odd and $l=1, 3, 5, \cdots, k-2$. Hence,
$$\left\{\aligned
  a_1(z) &=q^{\frac{121}{104}} (q; q^{13})(q^{12}; q^{13})(q^{13}; q^{13}),\\
  a_2(z) &=q^{\frac{49}{104}} (q^3; q^{13})(q^{10}; q^{13})(q^{13}; q^{13}),\\
  a_3(z) &=q^{\frac{25}{104}} (q^9; q^{13})(q^4; q^{13})(q^{13}; q^{13}),\\
  a_4(z) &=-q^{\frac{9}{104}} (q^5; q^{13})(q^8; q^{13})(q^{13}; q^{13}),\\
  a_5(z) &=q^{\frac{81}{104}} (q^2; q^{13})(q^{11}; q^{13})(q^{13}; q^{13}),\\
  a_6(z) &=q^{\frac{1}{104}} (q^6; q^{13})(q^7; q^{13})(q^{13}; q^{13}).
\endaligned\right.\eqno{(2.12)}$$
Ramanujan's general theta function $f(a, b)$ is given by
$$f(a, b):=\sum_{n=-\infty}^{\infty} a^{n(n+1)/2} b^{n(n-1)/2}, \quad |ab|<1,$$
and
$$f(-q):=f(-q, -q^2)=\sum_{n=-\infty}^{\infty} (-1)^n q^{n(3n-1)/2}=(q; q)_{\infty}.$$
In his notebooks (see \cite{R5}, p.326, \cite{R6}, p.244 and
\cite{B1}, p.372), Ramanujan obtained his modular equations of
degree $13$.

{\bf Theorem 2.4}. {\it Define
$$\mu_1=\frac{f(-q^4, -q^9)}{q^{7/13} f(-q^2, -q^{11})}, \quad
  \mu_2=\frac{f(-q^6, -q^7)}{q^{6/13} f(-q^3, -q^{10})},$$
$$\mu_3=\frac{f(-q^2, -q^{11})}{q^{5/13} f(-q, -q^{12})}, \quad
  \mu_4=\frac{f(-q^5, -q^8)}{q^{2/13} f(-q^4, -q^9)},$$
$$\mu_5=\frac{q^{5/13} f(-q^3, -q^{10})}{f(-q^5, -q^8)}, \quad
  \mu_6=\frac{q^{15/13} f(-q, -q^{12})}{f(-q^6, -q^7)}.$$
Then
$$1+\frac{f^2(-q)}{q f^2(-q^{13})}=\mu_1 \mu_2-\mu_3 \mu_5-\mu_4 \mu_6,\eqno{(2.13)}$$
$$-4-\frac{f^2(-q)}{q f^2(-q^{13})}=\frac{1}{\mu_1 \mu_2}-\frac{1}{\mu_3 \mu_5}-\frac{1}{\mu_4 \mu_6},\eqno{(2.14)}$$
$$3+\frac{f^2(-q)}{q f^2(-q^{13})}=\mu_2 \mu_3 \mu_4-\mu_1 \mu_5 \mu_6,\eqno{(2.15)}$$
where $\mu_1 \mu_2 \mu_3 \mu_4 \mu_5 \mu_6=1$.}

  Let
$$G_m(z):=G_{m, p}(z):=(-1)^m q^{m(3m-p)/(2p^2)} \frac{f(-q^{2m/p},
          -q^{1-2m/p})}{f(-q^{m/p}, -q^{1-m/p})}$$
where $m$ is a positive integer, $p=13$ and $q=e^{2 \pi i z}$.
Then the above three formulas (2.13), (2.14) and (2.15) are
equivalent to the following (see \cite{Ev} or \cite{B1},
p.373-375)
$$1+G_1(z) G_5(z)+G_2(z) G_3(z)+G_4(z) G_6(z)=-\frac{\eta^2(z/13)}{\eta^2(z)},\eqno{(2.16)}$$
$$\frac{1}{G_1(z) G_5(z)}+\frac{1}{G_2(z) G_3(z)}+\frac{1}{G_4(z) G_6(z)}
 =4+\frac{\eta^2(z/13)}{\eta^2(z)},\eqno{(2.17)}$$
$$G_1(z) G_3(z) G_4(z)-\frac{1}{G_1(z) G_3(z) G_4(z)}
 =3+\frac{\eta^2(z/13)}{\eta^2(z)},\eqno{(2.18)}$$
where $G_1(z) G_2(z) G_3(z) G_4(z) G_5(z) G_6(z)=-1$. Moreover,
there is the following formula (see \cite{Ev} or \cite{B1},
p.375-376): for $t=q^{1/13}$,
$$\frac{1}{(t^2)_{\infty} (t^3)_{\infty} (t^{10})_{\infty} (t^{11})_{\infty}}
 +\frac{t}{(t^4)_{\infty} (t^6)_{\infty} (t^7)_{\infty} (t^9)_{\infty}}
 =\frac{1}{(t)_{\infty} (t^5)_{\infty} (t^8)_{\infty} (t^{12})_{\infty}},\eqno{(2.19)}$$
which is equivalent to, for $p=13$, $G_1^{-1}(z)
G_5^{-1}(z)+G_4(z) G_6(z)=1$.

  It is known that (see \cite{Ev})
$$G(m; z)=(-1)^m F(2m/p; z)/F(m/p; z),$$
where
$$\aligned
  F(u; z) &=-i \sum_{k=-\infty}^{\infty} (-1)^k q^{(k+u+1/2)^2/2}\\
 &=-i q^{(u+1/2)^2/2} \prod_{m=1}^{\infty} (1-q^{m+u})(1-q^{m-1-u})
  (1-q^m)
\endaligned$$
satisfies that $F(u+1; z)=-F(u; z)$ and $F(-u; z)=-F(u; z)$. We
have
$$G_1(z) G_5(z)=\frac{F(\frac{2}{13}; z) F(\frac{3}{13}; z)}
                     {F(\frac{1}{13}; z) F(\frac{5}{13}; z)},$$
where
$$F\left(\frac{1}{13}; z\right)=-i q^{(\frac{1}{13}+\frac{1}{2})^2/2}
  \prod_{m=1}^{\infty} (1-q^{m+\frac{1}{13}})(1-q^{m-1-\frac{1}{13}})(1-q^m),$$
$$F\left(\frac{5}{13}; z\right)=-i q^{(\frac{5}{13}+\frac{1}{2})^2/2}
  \prod_{m=1}^{\infty} (1-q^{m+\frac{5}{13}})(1-q^{m-1-\frac{5}{13}})(1-q^m),$$
$$F\left(\frac{2}{13}; z\right)=-i q^{(\frac{2}{13}+\frac{1}{2})^2/2}
  \prod_{m=1}^{\infty} (1-q^{m+\frac{2}{13}})(1-q^{m-1-\frac{2}{13}})(1-q^m),$$
$$F\left(\frac{3}{13}; z\right)=-i q^{(\frac{3}{13}+\frac{1}{2})^2/2}
  \prod_{m=1}^{\infty} (1-q^{m+\frac{3}{13}})(1-q^{m-1-\frac{3}{13}})(1-q^m).$$
Note that
$$\prod_{m=1}^{\infty} (1-q^{m+\frac{1}{13}})=\frac{1}{1-t} (t)_{\infty}, \quad
  \prod_{m=1}^{\infty} (1-q^{m-1-\frac{1}{13}})=\frac{t-1}{t} (t^{12})_{\infty},$$
$$\prod_{m=1}^{\infty} (1-q^{m+\frac{5}{13}})=\frac{1}{1-t^5} (t^5)_{\infty}, \quad
  \prod_{m=1}^{\infty} (1-q^{m-1-\frac{5}{13}})=\frac{t^5-1}{t^5} (t^8)_{\infty},$$
$$\prod_{m=1}^{\infty} (1-q^{m+\frac{2}{13}})=\frac{1}{1-t^2} (t^2)_{\infty}, \quad
  \prod_{m=1}^{\infty} (1-q^{m-1-\frac{2}{13}})=\frac{t^2-1}{t^2} (t^{11})_{\infty},$$
$$\prod_{m=1}^{\infty} (1-q^{m+\frac{3}{13}})=\frac{1}{1-t^3} (t^3)_{\infty}, \quad
  \prod_{m=1}^{\infty} (1-q^{m-1-\frac{3}{13}})=\frac{t^3-1}{t^3} (t^{10})_{\infty},$$
where $(x)_{\infty}=\prod_{m=0}^{\infty} (1-xq^m)$ and
$t=q^{\frac{1}{13}}$. Thus,
$$G_1(z) G_5(z)=\frac{(t^3)_{\infty} (t^{10})_{\infty} (t^2)_{\infty} (t^{11})_{\infty}}
                {(t)_{\infty} (t^{12})_{\infty} (t^5)_{\infty}(t^8)_{\infty}}.$$
Similarly,
$$G_2(z) G_3(z)=-\frac{(t^9)_{\infty} (t^4)_{\infty} (t^6)_{\infty} (t^7)_{\infty}}
                {t (t^3)_{\infty} (t^{10})_{\infty} (t^2)_{\infty} (t^{11})_{\infty}}, \quad
  G_4(z) G_6(z)=\frac{t (t)_{\infty} (t^{12})_{\infty} (t^5)_{\infty} (t^8)_{\infty}}
                {(t^9)_{\infty} (t^4)_{\infty} (t^6)_{\infty}(t^7)_{\infty}}.$$
Note that $(t^k)_{\infty}=(t^k; t^{13})$ for $1 \leq k \leq 13$.
The above formula (2.16) is equivalent to
$$\aligned
 &t (t^3; t^{13})^2 (t^{10}; t^{13})^2 (t^2; t^{13})^2 (t^{11}; t^{13})^2
  (t^9; t^{13}) (t^4; t^{13}) (t^6; t^{13}) (t^7; t^{13})\\
-&(t^9; t^{13})^2 (t^4; t^{13})^2 (t^6; t^{13})^2 (t^7; t^{13})^2
  (t; t^{13}) (t^{12}; t^{13}) (t^5; t^{13}) (t^8; t^{13})\\
+&t^2 (t; t^{13})^2 (t^{12}; t^{13})^2 (t^5; t^{13})^2 (t^8;
t^{13})^2
  (t^3; t^{13}) (t^{10}; t^{13}) (t^2; t^{13}) (t^{11}; t^{13})\\
=&[-1-\eta^2(z/13)/\eta^2(z)] t (t; t^{13}) (t^2; t^{13}) \cdots
  (t^{12}; t^{13}), \quad t=q^{1/13}=e^{\frac{2 \pi i z}{13}}.
\endaligned$$
On the other hand,
$$\aligned
  a_2(z)^2 a_5(z)^2 a_3(z) a_6(z)
=&q^{\frac{11}{4}} (q^3; q^{13})^2 (q^{10}; q^{13})^2 (q^2; q^{13})^2 (q^{11}; q^{13})^2\\
 &\times (q^9; q^{13}) (q^4; q^{13}) (q^6; q^{13}) (q^7; q^{13}) (q^{13}; q^{13})^6,
\endaligned$$
$$\aligned
  a_3(z)^2 a_6(z)^2 a_1(z) a_4(z)
=&-q^{\frac{7}{4}} (q^9; q^{13})^2 (q^4; q^{13})^2 (q^6; q^{13})^2 (q^7; q^{13})^2\\
 &\times (q; q^{13}) (q^{12}; q^{13}) (q^5; q^{13}) (q^8; q^{13}) (q^{13}; q^{13})^6,
\endaligned$$
$$\aligned
  a_1(z)^2 a_4(z)^2 a_2(z) a_5(z)
=&q^{\frac{15}{4}} (q; q^{13})^2 (q^{12}; q^{13})^2 (q^5; q^{13})^2 (q^8; q^{13})^2\\
 &\times (q^3; q^{13}) (q^{10}; q^{13}) (q^2; q^{13}) (q^{11}; q^{13}) (q^{13}; q^{13})^6,
\endaligned$$
This implies that,
$$\aligned
 &a_1(z)^2 a_4(z)^2 a_2(z) a_5(z)+a_2(z)^2 a_5(z)^2 a_3(z) a_6(z)
  +a_3(z)^2 a_6(z)^2 a_1(z) a_4(z)\\
=&[-1-\eta^2(z)/\eta^2(13 z)] \eta(z) \eta(13 z)^5.
\endaligned\eqno{(2.20)}$$
Similarly, the other three formulas (2.17), (2.18) and (2.19) are
equivalent to the following:
$$\aligned
 &a_1(z)^2 a_4(z)^2 a_3(z) a_6(z)+a_2(z)^2 a_5(z)^2 a_1(z) a_4(z)
  +a_3(z)^2 a_6(z)^2 a_2(z) a_5(z)\\
=&[4+\eta^2(z)/\eta^2(13 z)] \eta(z) \eta(13 z)^5.
\endaligned\eqno{(2.21)}$$
$$a_4(z)^2 a_5(z)^2 a_6(z)^2-a_1(z)^2 a_2(z)^2 a_3(z)^2
 =[3+\eta^2(z)/\eta^2(13 z)] \eta(z) \eta(13 z)^5,\eqno{(2.22)}$$
and
$$\frac{1}{a_1(z) a_4(z)}+\frac{1}{a_2(z) a_5(z)}+\frac{1}{a_3(z) a_6(z)}=0,\eqno{(2.23)}$$
where $a_1(z) a_2(z) a_3(z) a_4(z) a_5(z) a_6(z)=-\eta(z) \eta(13
z)^5$.

  Let ${\bold A}=(a_1(z), a_2(z), a_3(z), a_4(z), a_5(z), a_6(z))^{T}$.
The significance of our six dimensional representation of $PSL(2,
13)$ comes from the following:

{\bf Proposition 2.5}. {\it If $z \in {\Bbb H}$, then the
following relations hold:
$${\bold A}(z+1)=e^{-\frac{3 \pi i}{4}} T {\bold A}(z), \quad
  {\bold A}\left(-\frac{1}{z}\right)=e^{\frac{\pi i}{4}} \sqrt{z} S {\bold A}(z),\eqno{(2.24)}$$
where
$$T=\left(\matrix
  \zeta^7 &   &   &   &   &    \\
  & \zeta^{11} &  &   &   &    \\
  &    & \zeta^8  &   &   &    \\
  &    &     & \zeta^6    &    \\
  &    &     &  & \zeta^2 &    \\
  &    &     &        &   & \zeta^5
\endmatrix\right),$$
$$S=-\frac{1}{\sqrt{13}} \left(\matrix
  \zeta^{12}-\zeta & \zeta^{10}-\zeta^3 & \zeta^4-\zeta^9 & \zeta^5-\zeta^8 & \zeta^2-\zeta^{11} & \zeta^6-\zeta^7\\
  \zeta^{10}-\zeta^3 & \zeta^4-\zeta^9 & \zeta^{12}-\zeta & \zeta^2-\zeta^{11} & \zeta^6-\zeta^7 & \zeta^5-\zeta^8\\
  \zeta^4-\zeta^9 & \zeta^{12}-\zeta & \zeta^{10}-\zeta^3 & \zeta^6-\zeta^7 & \zeta^5-\zeta^8 & \zeta^2-\zeta^{11}\\
  \zeta^5-\zeta^8 & \zeta^2-\zeta^{11} & \zeta^6-\zeta^7 & \zeta-\zeta^{12} & \zeta^3-\zeta^{10} & \zeta^9-\zeta^4\\
  \zeta^2-\zeta^{11} & \zeta^6-\zeta^7 & \zeta^5-\zeta^8 & \zeta^3-\zeta^{10} & \zeta^9-\zeta^4 & \zeta-\zeta^{12}\\
  \zeta^6-\zeta^7 & \zeta^5-\zeta^8 & \zeta^2-\zeta^{11} & \zeta^9-\zeta^4 & \zeta-\zeta^{12} & \zeta^3-\zeta^{10}
\endmatrix\right),$$
and $0<\text{arg} \sqrt{z} \leq \pi/2$.}

{\it Proof}. Recall that the following transformation formulas for
theta constants (see \cite{FK}, p.216-217):
$$\theta \left[\matrix \frac{2l+1}{k}\\ 1 \endmatrix\right](0, k(z+1))
 =\exp \left(\frac{\pi i}{4k}\right) \exp \left(\frac{\pi i}{k} (l^2+l) \right)
  \theta \left[\matrix \frac{2l+1}{k}\\ 1 \endmatrix\right](0, kz)$$
and
$$\aligned
 &\theta \left[\matrix \frac{2l+1}{k}\\ 1 \endmatrix\right](0, k(-1/z))
 =\frac{\sqrt{z}}{\sqrt{ik}} \exp \left(\frac{\pi i}{2k}\right) \times\\
 &\sum_{j=0}^{\frac{k-3}{2}} \left[\exp\left(2 \pi i \frac{l(j+1)}{k}\right)+
  \exp \left(-\frac{\pi i}{k}\right) \exp\left(-2 \pi i \frac{(l+1)j}{k}\right)\right]
  \theta \left[\matrix \frac{2j+1}{k}\\ 1 \endmatrix\right](0, kz),
\endaligned$$
where $k$ is odd and $l=0, 1, \cdots, \frac{k-3}{2}$. In our case,
$k=13$. By the first formula, we have
$$l=5, e^{\frac{\pi i}{52}} e^{\frac{30 \pi i}{13}}=e^{-\frac{3 \pi i}{4}} \zeta^7; \quad
  l=3, e^{\frac{\pi i}{52}} e^{\frac{12 \pi i}{13}}=e^{-\frac{3 \pi i}{4}} \zeta^{11}; \quad
  l=2, e^{\frac{\pi i}{52}} e^{\frac{6 \pi i}{13}}=e^{-\frac{3 \pi i}{4}} \zeta^8.$$
$$l=1, e^{\frac{\pi i}{52}} e^{\frac{2 \pi i}{13}}=e^{-\frac{3 \pi i}{4}} \zeta^6; \quad
  l=4, e^{\frac{\pi i}{52}} e^{\frac{20 \pi i}{13}}=e^{-\frac{3 \pi i}{4}} \zeta^2; \quad
  l=0, e^{\frac{\pi i}{52}}=e^{-\frac{3 \pi i}{4}} \zeta^5.$$
This gives that ${\bold A}(z+1)=e^{-\frac{3 \pi i}{4}} T {\bold
A}(z)$.

  By the second formula, for $l=5$, we have
$$\theta \left[\matrix \frac{11}{13}\\ 1 \endmatrix\right](0, -13/z)
=\frac{\sqrt{z}}{\sqrt{13 i}} e^{\frac{\pi i}{26}} \sum_{j=0}^{5}
 \left[e^{\frac{2 \pi i \cdot 5(j+1)}{13}}+e^{-\frac{\pi i}{13}}
 e^{-\frac{2 \pi i \cdot 6j}{13}}\right]
 \theta \left[\matrix \frac{2j+1}{13}\\ 1 \endmatrix\right](0, 13z).$$
where
$$\frac{1}{\sqrt{13 i}} e^{\frac{\pi i}{26}}
 =\frac{1}{\sqrt{13}} e^{-\frac{\pi i}{4}} e^{\frac{\pi i}{26}}
 =\frac{-1}{\sqrt{13}} e^{\frac{3 \pi i}{4}} e^{\frac{\pi i}{26}}
 =\frac{-1}{\sqrt{13}} e^{\frac{\pi i}{4}} e^{\frac{7 \pi i}{13}}.$$
Hence,
$$\aligned
  &\theta \left[\matrix \frac{11}{13}\\ 1 \endmatrix\right](0, -13/z)\\
 =&\frac{-1}{\sqrt{13}} \sqrt{z} e^{\frac{\pi i}{4}} e^{\frac{7 \pi i}{13}}
   [(e^{\frac{10 \pi i}{13}}+e^{-\frac{\pi i}{13}})
   \theta \left[\matrix \frac{1}{13}\\ 1 \endmatrix\right](0, 13z)+(
   e^{\frac{20 \pi i}{13}}+e^{-\pi i})
   \theta \left[\matrix \frac{3}{13}\\ 1 \endmatrix\right](0, 13z)+\\
  &+(e^{\frac{30 \pi i}{13}}+e^{-\frac{25 \pi i}{13}})
   \theta \left[\matrix \frac{5}{13}\\ 1 \endmatrix\right](0, 13z)+
   (e^{\frac{40 \pi i}{13}}+e^{-\frac{37 \pi i}{13}})
   \theta \left[\matrix \frac{7}{13}\\ 1 \endmatrix\right](0, 13z)+\\
  &+(e^{\frac{50 \pi i}{13}}+e^{-\frac{49 \pi i}{13}})
   \theta \left[\matrix \frac{9}{13}\\ 1 \endmatrix\right](0, 13z)+
   (e^{\frac{60 \pi i}{13}}+e^{-\frac{61 \pi i}{13}})
   \theta \left[\matrix \frac{11}{13}\\ 1 \endmatrix\right](0, 13z)].
\endaligned$$
Here,
$$e^{\frac{7 \pi i}{13}}(e^{\frac{10 \pi i}{13}}+e^{-\frac{\pi i}{13}})
 =e^{\frac{5 \pi i}{13}} (\zeta^6-\zeta^7), \quad
  e^{\frac{7 \pi i}{13}}(e^{\frac{20 \pi i}{13}}+e^{-\pi i})
 =-e^{\frac{4 \pi i}{13}} (\zeta^5-\zeta^8),$$
$$e^{\frac{7 \pi i}{13}}(e^{\frac{30 \pi i}{13}}+e^{-\frac{25 \pi i}{13}})
 =e^{\frac{3 \pi i}{13}} (\zeta^4-\zeta^9), \quad
  e^{\frac{7 \pi i}{13}}(e^{\frac{40 \pi i}{13}}+e^{-\frac{37 \pi i}{13}})
 =e^{\frac{2 \pi i}{13}} (\zeta^{10}-\zeta^3),$$
$$e^{\frac{7 \pi i}{13}}(e^{\frac{50 \pi i}{13}}+e^{-\frac{49 \pi i}{13}})
 =e^{\frac{\pi i}{13}} (\zeta^2-\zeta^{11}), \quad
  e^{\frac{7 \pi i}{13}}(e^{\frac{60 \pi i}{13}}+e^{-\frac{61 \pi i}{13}})
 =\zeta^{12}-\zeta.$$
This implies that
$$\aligned
  a_1\left(-\frac{1}{z}\right)
=&\frac{-1}{\sqrt{13}} \sqrt{z} e^{\frac{\pi i}{4}}
  [(\zeta^{12}-\zeta) a_1(z)+(\zeta^{10}-\zeta^3) a_2(z)+(\zeta^4-\zeta^9) a_3(z)+\\
 &+(\zeta^5-\zeta^8) a_4(z)+(\zeta^2-\zeta^{11}) a_5(z)+(\zeta^6-\zeta^7) a_6(z)].
\endaligned$$
For $l=3, 2, 1, 4, 0$, we can prove the similar formulas. This
gives that ${\bold A}\left(-\frac{1}{z}\right)=e^{\frac{\pi i}{4}}
\sqrt{z} S {\bold A}(z)$. This completes the proof of Proposition
2.5.

\flushpar $\qquad \qquad \qquad \qquad \qquad \qquad \qquad \qquad
\qquad \qquad \qquad \qquad \qquad \qquad \qquad \qquad \qquad
\qquad \quad \boxed{}$

\vskip 0.5 cm

\centerline{\bf 3. Seven-dimensional representations of $PSL(2,
                13)$, exotic modular equation}
\centerline{\bf and geometry of modular curve $X(13)$}

\vskip 0.5 cm

  We will construct a seven-dimensional representation of $PSL(2,
13)$ which is deduced from our six-dimensional representation. Let
us study the action of $S T^{\nu}$ on the five dimensional
projective space ${\Bbb P}^5=\{(z_1, z_2, z_3, z_4, z_5, z_6)\}$,
where $\nu=0, 1, \cdots, 12$. Put
$$\alpha=\zeta+\zeta^{12}-\zeta^5-\zeta^8, \quad
   \beta=\zeta^3+\zeta^{10}-\zeta^2-\zeta^{11}, \quad
   \gamma=\zeta^9+\zeta^4-\zeta^6-\zeta^7.$$
We find that
$$\aligned
  &13 ST^{\nu}(z_1) \cdot ST^{\nu}(z_4)\\
=&\beta z_1 z_4+\gamma z_2 z_5+\alpha z_3 z_6+\\
 &+\gamma \zeta^{\nu} z_1^2+\alpha \zeta^{9 \nu} z_2^2+\beta \zeta^{3 \nu} z_3^2
  -\gamma \zeta^{12 \nu} z_4^2-\alpha \zeta^{4 \nu} z_5^2-\beta \zeta^{10 \nu} z_6^2+\\
 &+(\alpha-\beta) \zeta^{5 \nu} z_1 z_2+(\beta-\gamma) \zeta^{6 \nu} z_2 z_3
  +(\gamma-\alpha) \zeta^{2 \nu} z_1 z_3+\\
 &+(\beta-\alpha) \zeta^{8 \nu} z_4 z_5+(\gamma-\beta) \zeta^{7 \nu} z_5 z_6
  +(\alpha-\gamma) \zeta^{11 \nu} z_4 z_6+\\
 &-(\alpha+\beta) \zeta^{\nu} z_3 z_4-(\beta+\gamma) \zeta^{9 \nu} z_1 z_5
  -(\gamma+\alpha) \zeta^{3 \nu} z_2 z_6+\\
 &-(\alpha+\beta) \zeta^{12 \nu} z_1 z_6-(\beta+\gamma) \zeta^{4 \nu} z_2 z_4
  -(\gamma+\alpha) \zeta^{10 \nu} z_3 z_5.
\endaligned$$
$$\aligned
  &13 ST^{\nu}(z_2) \cdot ST^{\nu}(z_5)\\
=&\gamma z_1 z_4+\alpha z_2 z_5+\beta z_3 z_6+\\
 &+\alpha \zeta^{\nu} z_1^2+\beta \zeta^{9 \nu} z_2^2+\gamma \zeta^{3 \nu} z_3^2
  -\alpha \zeta^{12 \nu} z_4^2-\beta \zeta^{4 \nu} z_5^2-\gamma \zeta^{10 \nu} z_6^2+\\
 &+(\beta-\gamma) \zeta^{5 \nu} z_1 z_2+(\gamma-\alpha) \zeta^{6 \nu} z_2 z_3
  +(\alpha-\beta) \zeta^{2 \nu} z_1 z_3+\\
 &+(\gamma-\beta) \zeta^{8 \nu} z_4 z_5+(\alpha-\gamma) \zeta^{7 \nu} z_5 z_6
  +(\beta-\alpha) \zeta^{11 \nu} z_4 z_6+\\
 &-(\beta+\gamma) \zeta^{\nu} z_3 z_4-(\gamma+\alpha) \zeta^{9 \nu} z_1 z_5
  -(\alpha+\beta) \zeta^{3 \nu} z_2 z_6+\\
 &-(\beta+\gamma) \zeta^{12 \nu} z_1 z_6-(\gamma+\alpha) \zeta^{4 \nu} z_2 z_4
  -(\alpha+\beta) \zeta^{10 \nu} z_3 z_5.
\endaligned$$
$$\aligned
  &13 ST^{\nu}(z_3) \cdot ST^{\nu}(z_6)\\
=&\alpha z_1 z_4+\beta z_2 z_5+\gamma z_3 z_6+\\
 &+\beta \zeta^{\nu} z_1^2+\gamma \zeta^{9 \nu} z_2^2+\alpha \zeta^{3 \nu} z_3^2
  -\beta \zeta^{12 \nu} z_4^2-\gamma \zeta^{4 \nu} z_5^2-\alpha \zeta^{10 \nu} z_6^2+\\
 &+(\gamma-\alpha) \zeta^{5 \nu} z_1 z_2+(\alpha-\beta) \zeta^{6 \nu} z_2 z_3
  +(\beta-\gamma) \zeta^{2 \nu} z_1 z_3+\\
 &+(\alpha-\gamma) \zeta^{8 \nu} z_4 z_5+(\beta-\alpha) \zeta^{7 \nu} z_5 z_6
  +(\gamma-\beta) \zeta^{11 \nu} z_4 z_6+\\
 &-(\gamma+\alpha) \zeta^{\nu} z_3 z_4-(\alpha+\beta) \zeta^{9 \nu} z_1 z_5
  -(\beta+\gamma) \zeta^{3 \nu} z_2 z_6+\\
 &-(\gamma+\alpha) \zeta^{12 \nu} z_1 z_6-(\alpha+\beta) \zeta^{4 \nu} z_2 z_4
  -(\beta+\gamma) \zeta^{10 \nu} z_3 z_5.
\endaligned$$
Note that $\alpha+\beta+\gamma=\sqrt{13}$, we find that
$$\aligned
  &\sqrt{13} \left[ST^{\nu}(z_1) \cdot ST^{\nu}(z_4)+ST^{\nu}(z_2) \cdot ST^{\nu}(z_5)+ST^{\nu}(z_3) \cdot ST^{\nu}(z_6)\right]\\
 =&(z_1 z_4+z_2 z_5+z_3 z_6)+(\zeta^{\nu} z_1^2+\zeta^{9 \nu} z_2^2+\zeta^{3 \nu} z_3^2)
  -(\zeta^{12 \nu} z_4^2+\zeta^{4 \nu} z_5^2+\zeta^{10 \nu} z_6^2)+\\
  &-2(\zeta^{\nu} z_3 z_4+\zeta^{9 \nu} z_1 z_5+\zeta^{3 \nu} z_2 z_6)
   -2(\zeta^{12 \nu} z_1 z_6+\zeta^{4 \nu} z_2 z_4+\zeta^{10 \nu} z_3 z_5).
\endaligned$$
Let
$$\varphi_{\infty}(z_1, z_2, z_3, z_4, z_5, z_6)=\sqrt{13} (z_1 z_4+z_2 z_5+z_3 z_6)\eqno{(3.1)}$$
and
$$\varphi_{\nu}(z_1, z_2, z_3, z_4, z_5, z_6)=\varphi_{\infty}(ST^{\nu}(z_1, z_2, z_3, z_4, z_5, z_6))\eqno{(3.2)}$$
for $\nu=0, 1, \cdots, 12$. Then
$$\aligned
  \varphi_{\nu}
=&(z_1 z_4+z_2 z_5+z_3 z_6)+\zeta^{\nu} (z_1^2-2 z_3 z_4)+\zeta^{4 \nu} (-z_5^2-2 z_2 z_4)+\\
 &+\zeta^{9 \nu} (z_2^2-2 z_1 z_5)+\zeta^{3 \nu} (z_3^2-2 z_2 z_6)+
   \zeta^{12 \nu} (-z_4^2-2 z_1 z_6)+\zeta^{10 \nu} (-z_6^2-2 z_3 z_5).
\endaligned\eqno{(3.3)}$$
This leads us to define the following senary quadratic forms
(quadratic forms in six variables):
$$\left\{\aligned
  {\Bbb A}_0 &=z_1 z_4+z_2 z_5+z_3 z_6,\\
  {\Bbb A}_1 &=z_1^2-2 z_3 z_4,\\
  {\Bbb A}_2 &=-z_5^2-2 z_2 z_4,\\
  {\Bbb A}_3 &=z_2^2-2 z_1 z_5,\\
  {\Bbb A}_4 &=z_3^2-2 z_2 z_6,\\
  {\Bbb A}_5 &=-z_4^2-2 z_1 z_6,\\
  {\Bbb A}_6 &=-z_6^2-2 z_3 z_5.
\endaligned\right.\eqno{(3.4)}$$
Hence,
$$\sqrt{13} ST^{\nu}({\Bbb A}_0)={\Bbb A}_0+\zeta^{\nu} {\Bbb A}_1+\zeta^{4 \nu} {\Bbb A}_2+
  \zeta^{9 \nu} {\Bbb A}_3+\zeta^{3 \nu} {\Bbb A}_4+\zeta^{12 \nu} {\Bbb A}_5+\zeta^{10 \nu} {\Bbb A}_6.\eqno{(3.5)}$$
Let
$$\left\{\aligned
  p_1 &=\zeta^2+\zeta^{11}-2+2(\zeta+\zeta^{12}-\zeta^9-\zeta^4),\\
  p_2 &=2-\zeta^9-\zeta^4+2(\zeta^5+\zeta^8-\zeta^2-\zeta^{11}),\\
  p_3 &=\zeta^6+\zeta^7-2+2(\zeta^3+\zeta^{10}-\zeta-\zeta^{12}),\\
  p_4 &=\zeta^5+\zeta^8-2+2(\zeta^9+\zeta^4-\zeta^3-\zeta^{10}),\\
  p_5 &=2-\zeta^3-\zeta^{10}+2(\zeta^6+\zeta^7-\zeta^5-\zeta^8),\\
  p_6 &=2-\zeta-\zeta^{12}+2(\zeta^2+\zeta^{11}-\zeta^6-\zeta^7).
\endaligned\right.\eqno{(3.6)}$$
We find that
$$\left\{\aligned
  13 S({\Bbb A}_1) &=2 \sqrt{13} {\Bbb A}_0+p_1 {\Bbb A}_1+p_2 {\Bbb A}_2+
                     p_3 {\Bbb A}_3+p_4 {\Bbb A}_4+p_5 {\Bbb A}_5+p_6 {\Bbb A}_6,\\
  13 S({\Bbb A}_2) &=2 \sqrt{13} {\Bbb A}_0+p_2 {\Bbb A}_1+p_4 {\Bbb A}_2+
                     p_6 {\Bbb A}_3+p_5 {\Bbb A}_4+p_3 {\Bbb A}_5+p_1 {\Bbb A}_6,\\
  13 S({\Bbb A}_3) &=2 \sqrt{13} {\Bbb A}_0+p_3 {\Bbb A}_1+p_6 {\Bbb A}_2+
                     p_4 {\Bbb A}_3+p_1 {\Bbb A}_4+p_2 {\Bbb A}_5+p_5 {\Bbb A}_6,\\
  13 S({\Bbb A}_4) &=2 \sqrt{13} {\Bbb A}_0+p_4 {\Bbb A}_1+p_5 {\Bbb A}_2+
                     p_1 {\Bbb A}_3+p_3 {\Bbb A}_4+p_6 {\Bbb A}_5+p_2 {\Bbb A}_6,\\
  13 S({\Bbb A}_5) &=2 \sqrt{13} {\Bbb A}_0+p_5 {\Bbb A}_1+p_3 {\Bbb A}_2+
                     p_2 {\Bbb A}_3+p_6 {\Bbb A}_4+p_1 {\Bbb A}_5+p_4 {\Bbb A}_6,\\
  13 S({\Bbb A}_6) &=2 \sqrt{13} {\Bbb A}_0+p_6 {\Bbb A}_1+p_1 {\Bbb A}_2+
                     p_5 {\Bbb A}_3+p_2 {\Bbb A}_4+p_4 {\Bbb A}_5+p_3 {\Bbb A}_6,\\
\endaligned\right.\eqno{(3.7)}$$
Note that
$$\left\{\aligned
  p_1 &=\sqrt{13} (\zeta^2+\zeta^{11}), \\
  p_2 &=\sqrt{13} (\zeta^9+\zeta^4),\\
  p_3 &=\sqrt{13} (\zeta^6+\zeta^7),\\
  p_4 &=\sqrt{13} (\zeta^5+\zeta^8),\\
  p_5 &=\sqrt{13} (\zeta^3+\zeta^{10}),\\
  p_6 &=\sqrt{13} (\zeta+\zeta^{12}).
\endaligned\right.\eqno{(3.8)}$$
We obtain a seven-dimensional representation of the simple group
$PSL(2, 13) \cong \langle \widetilde{S}, \widetilde{T} \rangle$
which is deduced from the actions of $S$ and $T$ on the basis
$({\Bbb A}_0, {\Bbb A}_1, {\Bbb A}_2, {\Bbb A}_3, {\Bbb A}_4,
{\Bbb A}_5, {\Bbb A}_6)$. Here
$$\widetilde{S}=\frac{1}{\sqrt{13}} \left(\matrix
  1 & 1                  & 1               & 1               & 1               & 1                  & 1\\
  2 & \zeta^2+\zeta^{11} & \zeta^9+\zeta^4 & \zeta^6+\zeta^7 & \zeta^5+\zeta^8 & \zeta^3+\zeta^{10} & \zeta+\zeta^{12}\\
  2 & \zeta^9+\zeta^4 & \zeta^5+\zeta^8 & \zeta+\zeta^{12} & \zeta^3+\zeta^{10} & \zeta^6+\zeta^7 & \zeta^2+\zeta^{11}\\
  2 & \zeta^6+\zeta^7 & \zeta+\zeta^{12} & \zeta^5+\zeta^8 & \zeta^2+\zeta^{11} & \zeta^9+\zeta^4 & \zeta^3+\zeta^{10}\\
  2 & \zeta^5+\zeta^8 & \zeta^3+\zeta^{10} & \zeta^2+\zeta^{11} & \zeta^6+\zeta^7 & \zeta+\zeta^{12} & \zeta^9+\zeta^4\\
  2 & \zeta^3+\zeta^{10} & \zeta^6+\zeta^7 & \zeta^9+\zeta^4 & \zeta+\zeta^{12} & \zeta^2+\zeta^{11} & \zeta^5+\zeta^8\\
  2 & \zeta+\zeta^{12} & \zeta^2+\zeta^{11} & \zeta^3+\zeta^{10} & \zeta^9+\zeta^4 & \zeta^5+\zeta^8 & \zeta^6+\zeta^7
\endmatrix\right),\eqno{(3.9)}$$
and
$$\widetilde{T}=\left(\matrix
  1 &       &         &         &         &            &           \\
    & \zeta &         &         &         &            &           \\
    &       & \zeta^4 &         &         &            &           \\
    &       &         & \zeta^9 &         &            &           \\
    &       &         &         & \zeta^3 &            &           \\
    &       &         &         &         & \zeta^{12} &           \\
    &       &         &         &         &            & \zeta^{10}
\endmatrix\right).\eqno{(3.10)}$$

  Now, let us reproduce from \cite{CC} some of the character table of $G$.

$$\text{Table $1$. Some of the character table of $PSL(2, 13)$}$$
$$\matrix
            & 1A & 2A & 3A & 6A & 7A & 7B & 7C & 13A & 13B\\
     \chi_1 &  1 &  1 &  1 &  1 &  1 &  1 &  1 &   1 &   1\\
     \chi_2 &  7 & -1 &  1 & -1 &  0 &  0 &  0 & \frac{1-\sqrt{13}}{2} & \frac{1+\sqrt{13}}{2}\\
     \chi_3 &  7 & -1 &  1 & -1 &  0 &  0 &  0 & \frac{1+\sqrt{13}}{2} & \frac{1-\sqrt{13}}{2}\\
  \chi_{10} &  6 &  0 &  0 &  0 & -1 & -1 & -1 & \frac{-1+\sqrt{13}}{2} & \frac{-1-\sqrt{13}}{2}\\
  \chi_{11} &  6 &  0 &  0 &  0 & -1 & -1 & -1 & \frac{-1-\sqrt{13}}{2} & \frac{-1+\sqrt{13}}{2}\\
  \chi_{15} & 14 &  0 &  2 &  0 &  0 &  0 &  0 &   1 &   1
\endmatrix$$

  In \cite{K1}, Klein obtained the modular equation of degree fourteen,
which corresponds to the transformation of order thirteen:
$$\aligned
J: J-1: 1=&(\tau^2+5 \tau+13)(\tau^4+7 \tau^3+20 \tau^2+19 \tau+1)^3\\
         :&(\tau^2+6 \tau+13)(\tau^6+10 \tau^5+46 \tau^4+108 \tau^3+122 \tau^2+38 \tau-1)^2\\
         :&1728 \tau,
\endaligned$$
where the Hauptmodul
$\tau=(\eta(z)/\eta(13 z))^2$. Note that
$$\aligned
 &\tau^4+7 \tau^3+20 \tau^2+19 \tau+1\\
=&\left(\tau^2+\frac{7+\sqrt{13}}{2} \tau+\frac{11+3
  \sqrt{13}}{2}\right)\left(\tau^2+\frac{7-\sqrt{13}}{2} \tau+\frac{11-3 \sqrt{13}}{2}\right),
\endaligned$$
$$\aligned
 &\tau^6+10 \tau^5+46 \tau^4+108 \tau^3+122 \tau^2+38 \tau-1\\
=&\left(\tau^3+5 \tau^2+\frac{21-\sqrt{13}}{2}
  \tau+\frac{3+\sqrt{13}}{2}\right) \left(\tau^3+5 \tau^2+\frac{21+\sqrt{13}}{2} \tau
  +\frac{3-\sqrt{13}}{2}\right).
\endaligned$$
In \cite{De1} and \cite{De2}, Deligne gave a modular
interpretation of why such factorizations exist. Let us confine
our thought to an especially important result which Jacobi had
established as early as $1829$ in his ``Notices sur les fonctions
elliptiques'' (see \cite{K}). Jacobi there considered, instead of
the modular equation, the so-called multiplier-equation, together
with other equations equivalent to it, and found that their
$(n+1)$ roots are composed in a simple manner of $\frac{n+1}{2}$
elements, with the help of merely numerical irrationalities.
Namely, if we denote these elements by ${\Bbb A}_0$, ${\Bbb A}_1$,
$\cdots$, ${\Bbb A}_{\frac{n-1}{2}}$, and further, for the roots
$z$ of the equation under consideration, apply the indices
employed by Galois, we have, with appropriate determination of the
square root occurring on the left-hand side:
$$\left\{\aligned
  \sqrt{z_{\infty}} &=\sqrt{(-1)^{\frac{n-1}{2}} \cdot n} \cdot {\Bbb A}_0,\\
  \sqrt{z_{\nu}} &={\Bbb A}_0+\epsilon^{\nu} {\Bbb A}_1+\epsilon^{4 \nu} {\Bbb A}_2
  +\cdots+\epsilon^{(\frac{n-1}{2})^2 \nu} {\Bbb A}_{\frac{n-1}{2}}
\endaligned\right.$$
for $\nu=0, 1, \cdots, n-1$ and $\epsilon=e^{\frac{2 \pi i}{n}}$.
Jacobi had himself emphasized the special significance of his
result by adding to his short communication: ``C'est un
th\'{e}or\`{e}me des plus importants dans la th\'{e}orie
alg\'{e}brique de la transformation et de la division des
fonctions elliptiques.'' Now, we give the Jacobian equation of
degree fourteen which corresponds to the above modular equation of
degree fourteen.

  Let $H:=Q^5 P^2 \cdot P^2 Q^6 P^8 \cdot Q^5 P^2 \cdot P^3 Q$ where
$$Q^5 P^2=-\frac{1}{\sqrt{13}} \left(\matrix
  \zeta^8-\zeta^5 & \zeta^4-\zeta^8 & \zeta^2-\zeta & \zeta^4-\zeta^6 & \zeta^9-\zeta^2 & \zeta-\zeta^6\\
  \zeta^5-\zeta^9 & \zeta^7-\zeta^6 & \zeta^{10}-\zeta^7 & \zeta^9-\zeta^2 & \zeta^{10}-\zeta^2 & \zeta^3-\zeta^5\\
  \zeta^{12}-\zeta^{11} & \zeta^6-\zeta^3 & \zeta^{11}-\zeta^2 & \zeta-\zeta^6 & \zeta^3-\zeta^5 & \zeta^{12}-\zeta^5\\
  \zeta^7-\zeta^9 & \zeta^{11}-\zeta^4 & \zeta^7-\zeta^{12} & \zeta^5-\zeta^8 & \zeta^9-\zeta^5 & \zeta^{11}-\zeta^{12}\\
  \zeta^{11}-\zeta^4 & \zeta^{11}-\zeta^3 & \zeta^8-\zeta^{10} & \zeta^8-\zeta^4 & \zeta^6-\zeta^7 & \zeta^3-\zeta^6\\
  \zeta^7-\zeta^{12} & \zeta^8-\zeta^{10} & \zeta^8-\zeta & \zeta-\zeta^2 & \zeta^7-\zeta^{10} & \zeta^2-\zeta^{11}
\endmatrix\right),$$
$$P^2 Q^6 P^8=-\frac{1}{\sqrt{13}} \left(\matrix
  \zeta^8-\zeta^5 & \zeta^{12}-\zeta^3 & \zeta^4-\zeta^3 & \zeta^2-\zeta^4 & \zeta^{12}-\zeta^5 & \zeta^{10}-\zeta^2\\
  \zeta^{10}-\zeta & \zeta^7-\zeta^6 & \zeta^4-\zeta & \zeta^{12}-\zeta^5 & \zeta^5-\zeta^{10} & \zeta^4-\zeta^6\\
  \zeta^{10}-\zeta^9 & \zeta^{12}-\zeta^9 & \zeta^{11}-\zeta^2 & \zeta^{10}-\zeta^2 & \zeta^4-\zeta^6 & \zeta^6-\zeta^{12}\\
  \zeta^9-\zeta^{11} & \zeta^8-\zeta & \zeta^{11}-\zeta^3 & \zeta^5-\zeta^8 & \zeta-\zeta^{10} & \zeta^9-\zeta^{10}\\
  \zeta^8-\zeta & \zeta^3-\zeta^8 & \zeta^7-\zeta^9 & \zeta^3-\zeta^{12} & \zeta^6-\zeta^7 & \zeta^9-\zeta^{12}\\
  \zeta^{11}-\zeta^3 & \zeta^7-\zeta^9 & \zeta-\zeta^7 & \zeta^3-\zeta^4 & \zeta-\zeta^4 & \zeta^2-\zeta^{11}
\endmatrix\right),$$
and
$$\aligned
 &Q^5 P^2 \cdot P^2 Q^6 P^8 \cdot Q^5 P^2\\
=&-\frac{1}{\sqrt{13}} \left(\matrix
  \zeta^7-\zeta^6 & \zeta^8-\zeta^5 & \zeta^{11}-\zeta^2 & \zeta^4-\zeta^9 & \zeta^{12}-\zeta & \zeta^{10}-\zeta^3\\
  \zeta^8-\zeta^5 & \zeta^{11}-\zeta^2 & \zeta^7-\zeta^6 & \zeta^{12}-\zeta & \zeta^{10}-\zeta^3 & \zeta^4-\zeta^9\\
  \zeta^{11}-\zeta^2 & \zeta^7-\zeta^6 & \zeta^8-\zeta^5 & \zeta^{10}-\zeta^3 & \zeta^4-\zeta^9 & \zeta^{12}-\zeta\\
  \zeta^4-\zeta^9 & \zeta^{12}-\zeta & \zeta^{10}-\zeta^3 & \zeta^6-\zeta^7 & \zeta^5-\zeta^8 & \zeta^2-\zeta^{11}\\
  \zeta^{12}-\zeta & \zeta^{10}-\zeta^3 & \zeta^4-\zeta^9 & \zeta^5-\zeta^8 & \zeta^2-\zeta^{11} & \zeta^6-\zeta^7\\
  \zeta^{10}-\zeta^3 & \zeta^4-\zeta^9 & \zeta^{12}-\zeta & \zeta^2-\zeta^{11} & \zeta^6-\zeta^7 & \zeta^5-\zeta^8
\endmatrix\right).\endaligned$$
Then
$$H=\left(\matrix
  0 &  0 &  0 & 0 & 0 & 1\\
  0 &  0 &  0 & 1 & 0 & 0\\
  0 &  0 &  0 & 0 & 1 & 0\\
  0 &  0 & -1 & 0 & 0 & 0\\
 -1 &  0 &  0 & 0 & 0 & 0\\
  0 & -1 &  0 & 0 & 0 & 0
\endmatrix\right).\eqno{(3.11)}$$
Note that
$$H^2=\left(\matrix
  0 & -1 &  0 &  0 &  0 &  0\\
  0 &  0 & -1 &  0 &  0 &  0\\
 -1 &  0 &  0 &  0 &  0 &  0\\
  0 &  0 &  0 &  0 & -1 &  0\\
  0 &  0 &  0 &  0 &  0 & -1\\
  0 &  0 &  0 & -1 &  0 &  0
\endmatrix\right), \quad
  H^3=\left(\matrix
  0 & 0 & 0 & -1 &  0 &  0\\
  0 & 0 & 0 &  0 & -1 &  0\\
  0 & 0 & 0 &  0 &  0 & -1\\
  1 & 0 & 0 &  0 &  0 &  0\\
  0 & 1 & 0 &  0 &  0 &  0\\
  0 & 0 & 1 &  0 &  0 &  0
\endmatrix\right)$$
and $H^6=-I$. In the projective coordinates, this means that
$H^6=1$. We have
$$H^{-1} T H=-T^4.$$
Thus, $\langle H, T \rangle \cong {\Bbb Z}_{13} \rtimes {\Bbb
Z}_6$. Hence, it is a maximal subgroup of order $78$ of $G$ with
index $14$ (see \cite{CC}). We find that $\varphi_{\infty}^2$ is
invariant under the action of the maximal subgroup $\langle H, T
\rangle$. Note that
$$\varphi_{\infty}=\sqrt{13} {\Bbb A}_0, \quad
  \varphi_{\nu}={\Bbb A}_0+\zeta^{\nu} {\Bbb A}_1+\zeta^{4 \nu} {\Bbb A}_2+
  \zeta^{9 \nu} {\Bbb A}_3+\zeta^{3 \nu} {\Bbb A}_4+\zeta^{12 \nu} {\Bbb A}_5+\zeta^{10 \nu} {\Bbb A}_6$$
for $\nu=0, 1, \cdots, 12$. Let $w=\varphi^2$,
$w_{\infty}=\varphi_{\infty}^2$ and $w_{\nu}=\varphi_{\nu}^2$.
Then $w_{\infty}$, $w_{\nu}$ for $\nu=0, \cdots, 12$ form an
algebraic equation of degree fourteen, which is just the Jacobian
equation of degree fourteen, whose roots are these $w_{\nu}$ and
$w_{\infty}$:
$$w^{14}+a_1 w^{13}+\cdots+ a_{13} w+a_{14}=0.$$
In particular, the coefficients
$$a_{14}=\varphi_{\infty}^2 \cdot \prod_{\nu=0}^{12} \varphi_{\nu}^2
 =13 {\Bbb A}_0^2 \prod_{\nu=0}^{12} ({\Bbb A}_0+\zeta^{\nu} {\Bbb A}_1+\zeta^{4 \nu} {\Bbb A}_2+
  \zeta^{9 \nu} {\Bbb A}_3+\zeta^{3 \nu} {\Bbb A}_4+\zeta^{12 \nu} {\Bbb A}_5+\zeta^{10 \nu} {\Bbb A}_6)^2,$$
and
$$-a_1=w_{\infty}+\sum_{\nu=0}^{12} w_{\nu}=
  26 ({\Bbb A}_0^2+{\Bbb A}_1 {\Bbb A}_5+{\Bbb A}_2 {\Bbb A}_3+{\Bbb A}_4 {\Bbb A}_6).\eqno{(3.12)}$$

  Recall that the modular equation of degree fourteen is also given by (see \cite{Hu})
$$\aligned
  z^{14}
 &+13 [2 \Delta^{\frac{1}{12}} z^{13}+25 \Delta^{\frac{2}{12}} z^{12}
  +196 \Delta^{\frac{3}{12}} z^{11}+1064 \Delta^{\frac{4}{12}} z^{10}\\
 &+4180 \Delta^{\frac{5}{12}} z^{9}+12086 \Delta^{\frac{6}{12}} z^{8}
  +25660 \Delta^{\frac{7}{12}} z^{7}+39182 \Delta^{\frac{8}{12}} z^{6}\\
 &+41140 \Delta^{\frac{9}{12}} z^{5}+27272 \Delta^{\frac{10}{12}} z^{4}
  +9604 \Delta^{\frac{11}{12}} z^{3}+1165 \Delta z^{2}]\\
 &+[746 \Delta-(12 g_2)^3] \Delta^{\frac{1}{12}} z+13 \Delta^{\frac{14}{12}}=0.
\endaligned$$
The coefficients of $z^{13}$ and the constant are just $26$ and
$13$.

  In fact, the invariant quadric
$$\Psi_2:={\Bbb A}_0^2+{\Bbb A}_1 {\Bbb A}_5+{\Bbb A}_2 {\Bbb A}_3+{\Bbb A}_4 {\Bbb A}_6
    =2 \Phi_4(z_1, z_2, z_3, z_4, z_5, z_6),\eqno{(3.13)}$$
where
$$\Phi_4:=(z_3 z_4^3+z_1 z_5^3+z_2 z_6^3)-(z_6 z_1^3+z_4 z_2^3+z_5 z_3^3)+
  3(z_1 z_2 z_4 z_5+z_2 z_3 z_5 z_6+z_3 z_1 z_6 z_4),\eqno{(3.14)}$$
Hence, the variety $\Psi_2=0$ is a quartic four-fold, which is
invariant under the action of the simple group $G$.

  Recall that the principal congruence subgroup of level $13$ is the
normal subgroup $\Gamma(13)$ of $\Gamma=PSL(2, {\Bbb Z})$ defined
by the exact sequence
$$1 \rightarrow \Gamma(13) \rightarrow \Gamma(1) @> f >> G \rightarrow 1\eqno{(3.15)}$$
where $f(\gamma) \equiv \gamma$ (mod $13$) for $\gamma \in
\Gamma=\Gamma(1)$. Then there is a representation
$$\rho: \Gamma \rightarrow PGL(6, {\Bbb C})\eqno{(3.16)}$$
with kernel $\Gamma(13)$ and leaving $\Phi_4$ invariant. It is
defined as follows: if $t=\left(\matrix 1 & 1\\ 0 & 1
\endmatrix\right)$ and $s=\left(\matrix 0 & -1\\ 1 & 0 \endmatrix\right)$,
then $\rho(t)=T$ and $\rho(s)=S$. To see that such a
representation exists, note that $\Gamma$ is defined by the
presentation
$$\langle s, t; s^2=(st)^3=1 \rangle$$
satisfied by $s$ and $t$ and we have proved that $S$ and $T$
satisfy these relations. Moreover, we have proved that $G$ is
defined by the presentation
$$\langle S, T; S^2=T^{13}=(ST)^3=1, (Q^3 P^4)^3=1 \rangle.\eqno{(3.17)}$$
Let $p=s t^{-1} s$ and $q=s t^3$. Then
$$h:=q^5 p^2 \cdot p^2 q^6 p^8 \cdot q^5 p^2 \cdot p^3 q
    =\left(\matrix
     4, 428, 249 & -10, 547, 030\\
    -11, 594, 791 & 27, 616, 019
     \endmatrix\right)\eqno{(3.18)}$$
satisfies that $\rho(h)=H$.

{\bf Theorem 3.1}. {\it The invariant quartic four-fold
$\Phi_4(z_1, z_2, z_3, z_4, z_5, z_6)=0$ can be parametrized by
the theta constants associated with $\Gamma(13):$
$$\Phi_4(a_1(z), a_2(z), a_3(z), a_4(z), a_5(z), a_6(z))=0.\eqno{(3.19)}$$}

{\it Proof}. Let $y_i(z)=\eta^3(z) a_i(z)$ $(1 \leq i \leq 6)$.
Then
$$Y(z):=\left(\matrix y_1(z)\\ y_2(z)\\ y_3(z)\\ y_4(z)\\ y_5(z)\\
        y_6(z) \endmatrix\right)
       =\eta^3(z) \left(\matrix a_1(z)\\ a_2(z)\\ a_3(z)\\
        a_4(z)\\ a_5(z)\\ a_6(z) \endmatrix\right)
       =\eta^3(z) {\bold A}(z).$$
It is well-known that $\eta(z)$ satisfies the following
transformation formulas
$$\eta(z+1)=e^{\frac{\pi i}{12}} \eta(z), \quad \text{and} \quad
  \eta\left(-\frac{1}{z}\right)=e^{-\frac{\pi i}{4}} \sqrt{z} \eta(z).$$
By Proposition 2.5, we have
$$Y(z+1)=e^{-\frac{\pi i}{2}} \rho(t) Y(z) \quad \text{and} \quad
  Y\left(-\frac{1}{z}\right)=e^{-\frac{\pi i}{2}} z^2 \rho(s) Y(z).$$
Define $j(\gamma, z):=cz+d$ if $z \in {\Bbb H}$ and
$\gamma=\left(\matrix a & b\\ c & d \endmatrix\right) \in
\Gamma(1)$. Hence,
$$Y(\gamma(z))=v(\gamma) j(\gamma, z)^2 \rho(\gamma) Y(z)
  \quad \text{for $\gamma \in \Gamma(1)$},$$
where $v(\gamma)=\pm 1$ or $\pm i$. Since $\Gamma(13)=\text{ker}$
$\rho$, we have
$$Y(\gamma(z))=v(\gamma) j(\gamma, z)^2 Y(z) \quad \text{for $\gamma \in \Gamma(13)$}.$$
This means that the functions $y_1(z)$, $y_2(z)$, $y_3(z)$,
$y_4(z)$, $y_5(z)$, $y_6(z)$ are modular forms of weight two for
$\Gamma(13)$ with the same multiplier $v(\gamma)=\pm 1 $ or $\pm
i$. Thus,
$$\Phi_4(Y(\gamma(z)))=v(\gamma)^4 j(\gamma, z)^8 \Phi_4(Y(z))
 =j(\gamma, z)^8 \Phi_4(Y(z)) \quad \text{for $\gamma \in \Gamma(13)$}.$$
Moreover, for $\gamma \in \Gamma(1)$,
$$\aligned
  &\Phi_4(Y(\gamma(z)))=\Phi_4(v(\gamma) j(\gamma, z)^2 \rho(\gamma) Y(z))\\
 =&v(\gamma)^4 j(\gamma, z)^8 \Phi_4(\rho(\gamma) Y(z))=j(\gamma, z)^8 \Phi_4(\rho(\gamma) Y(z)).
\endaligned$$
Note that $\rho(\gamma) \in \langle \rho(s), \rho(t) \rangle=G$
and $\Phi_4$ is a $G$-invariant polynomial, we have
$$\Phi_4(Y(\gamma(z)))=j(\gamma, z)^8 \Phi_4(Y(z)), \quad \text{for $\gamma \in \Gamma(1)$}.$$
This implies that $\Phi_4(Y(z))$ is a modular form of weight eight
for the full modular group $\Gamma(1)$. Moreover, we will show
that it is a cusp form. We have
$$a_3(z) a_4(z)^3=-q^{\frac{1}{2}} \sum_{n \in {\Bbb Z}} (-1)^n
  q^{\frac{1}{2}(13n^2+5n)} \left[\sum_{n \in {\Bbb Z}} (-1)^n
  q^{\frac{1}{2}(13n^2+3n)}\right]^3,$$
$$a_1(z) a_5(z)^3=q^{\frac{7}{2}} \sum_{n \in {\Bbb Z}} (-1)^n
  q^{\frac{1}{2}(13n^2+11n)} \left[\sum_{n \in {\Bbb Z}} (-1)^n
  q^{\frac{1}{2}(13n^2+9n)}\right]^3,$$
$$a_2(z) a_6(z)^3=q^{\frac{1}{2}} \sum_{n \in {\Bbb Z}} (-1)^n
  q^{\frac{1}{2}(13n^2+7n)} \left[\sum_{n \in {\Bbb Z}} (-1)^n
  q^{\frac{1}{2}(13n^2+n)}\right]^3,$$
$$a_6(z) a_1(z)^3=q^{\frac{7}{2}} \sum_{n \in {\Bbb Z}} (-1)^n
  q^{\frac{1}{2}(13n^2+n)} \left[\sum_{n \in {\Bbb Z}} (-1)^n
  q^{\frac{1}{2}(13n^2+11n)}\right]^3,$$
$$a_4(z) a_2(z)^3=-q^{\frac{3}{2}} \sum_{n \in {\Bbb Z}} (-1)^n
  q^{\frac{1}{2}(13n^2+3n)} \left[\sum_{n \in {\Bbb Z}} (-1)^n
  q^{\frac{1}{2}(13n^2+7n)}\right]^3,$$
$$a_5(z) a_3(z)^3=q^{\frac{3}{2}} \sum_{n \in {\Bbb Z}} (-1)^n
  q^{\frac{1}{2}(13n^2+9n)} \left[\sum_{n \in {\Bbb Z}} (-1)^n
  q^{\frac{1}{2}(13n^2+5n)}\right]^3,$$
$$\aligned
  a_1(z) a_2(z) a_4(z) a_5(z)
=&-q^{\frac{5}{2}} \sum_{n \in {\Bbb Z}} (-1)^n q^{\frac{1}{2}(13n^2+11n)}
  \sum_{n \in {\Bbb Z}} (-1)^n q^{\frac{1}{2}(13n^2+7n)}\\
 &\times \sum_{n \in {\Bbb Z}} (-1)^n  q^{\frac{1}{2}(13n^2+3n)}
  \sum_{n \in {\Bbb Z}} (-1)^n q^{\frac{1}{2}(13n^2+9n)},
\endaligned$$
$$\aligned
  a_2(z) a_3(z) a_5(z) a_6(z)
=&q^{\frac{3}{2}} \sum_{n \in {\Bbb Z}} (-1)^n q^{\frac{1}{2}(13n^2+7n)}
  \sum_{n \in {\Bbb Z}} (-1)^n q^{\frac{1}{2}(13n^2+5n)}\\
 &\times \sum_{n \in {\Bbb Z}} (-1)^n q^{\frac{1}{2}(13n^2+9n)}
  \sum_{n \in {\Bbb Z}} (-1)^n q^{\frac{1}{2}(13n^2+n)},
\endaligned$$
and
$$\aligned
  a_3(z) a_1(z) a_6(z) a_4(z)
=&-q^{\frac{3}{2}} \sum_{n \in {\Bbb Z}} q^{\frac{1}{2}(13n^2+5n)}
  \sum_{n \in {\Bbb Z}} (-1)^n q^{\frac{1}{2}(13n^2+11n)}\\
 &\times \sum_{n \in {\Bbb Z}} (-1)^n q^{\frac{1}{2}(13n^2+n)}
  \sum_{n \in {\Bbb Z}} (-1)^n q^{\frac{1}{2}(13n^2+3n)}.
\endaligned$$
Hence,
$$\Phi_4(a_1(z), a_2(z), a_3(z), a_4(z), a_5(z), a_6(z))
 =q^{\frac{1}{2}} \sum_{n=0}^{\infty} a_n q^n, \quad \text{where $a_n \in {\Bbb Z}$}.$$
On the other hand, $\eta(z)^{12}=q^{\frac{1}{2}}
\prod_{n=1}^{\infty} (1-q^n)^{12}$. We have
$$\Phi_4(y_1(z), y_2(z), y_3(z), y_4(z), y_5(z), y_6(z))
 =q \sum_{n=0}^{\infty} a_n q^n \prod_{n=1}^{\infty} (1-q^n)^{12}$$
is a cusp form of weight $8$ for the full modular group
$\Gamma=PSL(2, {\Bbb Z})$, but the only such form is zero. (See
\cite{Se}, Chapter VII for the complete description of cusp forms
on $\Gamma(1)$). This completes the proof of Theorem 3.1.

\flushpar $\qquad \qquad \qquad \qquad \qquad \qquad \qquad \qquad
\qquad \qquad \qquad \qquad \qquad \qquad \qquad \qquad \qquad
\qquad \quad \boxed{}$

{\bf Corollary 3.2}. {\it The following invariant decomposition
formula $($exotic modular equation$)$ holds:
$$\frac{\Psi_2(a_1(z), a_2(z), a_3(z), a_4(z), a_5(z), a_6(z))}
  {{\Bbb A}_0(a_1(z), a_2(z), a_3(z), a_4(z), a_5(z), a_6(z))^2}=0,\eqno{(3.20)}$$
where the quadric
$$\Psi_2(z_1, z_2, z_3, z_4, z_5, z_6):=
 {\Bbb A}_0^2+{\Bbb A}_1 {\Bbb A}_5+{\Bbb A}_2 {\Bbb A}_3+{\Bbb A}_4 {\Bbb A}_6\eqno{(3.21)}$$
is an invariant associated to $PSL(2, 13)$ and ${\Bbb A}_0^2$ is
invariant under the action of the image of a Borel subgroup of
$PSL(2, 13)$ $($i.e. a maximal subgroup of order $78$ of $PSL(2,
13)$$)$.}

{\it Proof}. This comes from Theorem 3.1 by noting that (3.13).
Thus, we complete the proof of Theorem 1.3.

\flushpar $\qquad \qquad \qquad \qquad \qquad \qquad \qquad \qquad
\qquad \qquad \qquad \qquad \qquad \qquad \qquad \qquad \qquad
\qquad \quad \boxed{}$

  Now, we will prove that
$$\aligned
      &\prod_{j=1}^{6} {\Bbb A}_j(a_1(z), a_2(z), a_3(z), a_4(z), a_5(z), a_6(z))\\
 \neq &-{\Bbb A}_0(a_1(z), a_2(z), a_3(z), a_4(z), a_5(z), a_6(z))^6.
\endaligned\eqno{(3.22)}$$
We have
$${\Bbb A}_0(a_1(z), a_2(z), a_3(z), a_4(z), a_5(z), a_6(z))=q^{\frac{1}{4}} (1+O(q)),$$
$${\Bbb A}_1(a_1(z), a_2(z), a_3(z), a_4(z), a_5(z), a_6(z))=q^{\frac{34}{104}} (2+O(q)),$$
$${\Bbb A}_2(a_1(z), a_2(z), a_3(z), a_4(z), a_5(z), a_6(z))=q^{\frac{58}{104}} (2+O(q)),$$
$${\Bbb A}_3(a_1(z), a_2(z), a_3(z), a_4(z), a_5(z), a_6(z))=q^{\frac{98}{104}} (1+O(q)),$$
$${\Bbb A}_4(a_1(z), a_2(z), a_3(z), a_4(z), a_5(z), a_6(z))=q^{\frac{50}{104}} (-1+O(q)),$$
$${\Bbb A}_5(a_1(z), a_2(z), a_3(z), a_4(z), a_5(z), a_6(z))=q^{\frac{18}{104}} (-1+O(q)),$$
$${\Bbb A}_6(a_1(z), a_2(z), a_3(z), a_4(z), a_5(z), a_6(z))=q^{\frac{2}{104}} (-1+O(q)).$$
Thus,
$$\prod_{j=1}^{6} {\Bbb A}_j(a_1(z), a_2(z), a_3(z), a_4(z), a_5(z), a_6(z))=q^{\frac{5}{2}} (-4+O(q)).$$
On the other hand,
$$-{\Bbb A}_0(a_1(z), a_2(z), a_3(z), a_4(z), a_5(z), a_6(z))^6=q^{\frac{3}{2}} (-1+O(q)).$$
This gives the proof of (3.22).

  In fact, Theorem 3.1 also gives the geometry of the classical modular
curve $X(13)$:

{\bf Corollary 3.3}. {\it The coordinates $(a_1(z), a_2(z),
a_3(z), a_4(z), a_5(z), a_6(z))$ map $X(13)$ into the quartic
four-fold $\Phi_4(z_1, z_2, z_3, z_4, z_5, z_6)=0$ in ${\Bbb C}
{\Bbb P}^5$.}

{\it Proof}. This comes from Theorem 2.3 and Theorem 3.1.
\flushpar $\qquad \qquad \qquad \qquad \qquad \qquad \qquad \qquad
\qquad \qquad \qquad \qquad \qquad \qquad \qquad \qquad \qquad
\qquad \quad \boxed{}$

\vskip 0.5 cm

\centerline{\bf 4. Fourteen-dimensional representations of $PSL(2,
                   13)$ and}
\centerline{\bf invariant decomposition formula}

\vskip 0.5 cm

  It is known that the principal congruence subgroup $\Gamma(13)$ is a
normal subgroup of $\Gamma_0(13)$, and $\Gamma_0(13)$ is generated
by $\Gamma(13)$ and the group generated by $H$ and $T$. The
quotient group $\Gamma_0(13)/\Gamma(13)$ is of order $78$ and is
isomorphic to the subgroup $\langle H, T \rangle$. The natural
projection defines a Galois covering $X(13) \to X_0(13)$ with
Galois group $\langle H, T \rangle$. Recall that $(\eta(z)/\eta(13
z))^2$ is the Hauptmodul of $X_0(13)$, that is, the unique (up to
multiplication by a constant) univalent modular function for
$\Gamma_0(13)$ holomorphic on ${\Bbb H}$ with a simple pole at
$z=i \infty$ and a simple zero at $z=0$.

  We will construct a fourteen-dimensional representation of $PSL(2,
13)$ which is deduced from our six-dimensional representation. It
should be emphasized that our fourteen-dimensional representation
is not Weil representation. In contrast with this fact, both our
six-dimensional and seven-dimensional representations of $G$ are
Weil representations, i.e., $\frac{p-1}{2}$-dimensional and
$\frac{p+1}{2}$-dimensional representations of $PSL(2, p)$,
respectively. In fact, what Klein used in his papers \cite{K},
\cite{K1}, \cite{K2}, \cite{K3} and \cite{K4} are all Weil
representations. Hence, our method is completely different from
Klein's method.

  To construct our fourteen-dimensional representation which is
generated under the action of $PSL$ by a specific vector in
$\text{Sym}^3$(six-dimensional representation), we begin with a
cubic polynomial $z_1 z_2 z_3$ and study the action of $ST^{\nu}$
($\nu$ mod $13$) on it. We have
$$\aligned
  &-13 \sqrt{13} ST^{\nu}(z_1) \cdot ST^{\nu}(z_2) \cdot ST^{\nu}(z_3)\\
 =&-\sqrt{\frac{-13-3 \sqrt{13}}{2}} (\zeta^{8 \nu} z_1^3+\zeta^{7 \nu} z_2^3+\zeta^{11 \nu} z_3^3)
   -\sqrt{\frac{-13+3 \sqrt{13}}{2}} (\zeta^{5 \nu} z_4^3+\zeta^{6 \nu} z_5^3+\zeta^{2 \nu} z_6^3)\\
  &-\sqrt{-13+2 \sqrt{13}} (\zeta^{12 \nu} z_1^2 z_2+\zeta^{4 \nu} z_2^2 z_3+\zeta^{10 \nu} z_3^2 z_1)\\
  &-\sqrt{-13-2 \sqrt{13}} (\zeta^{\nu} z_4^2 z_5+\zeta^{9 \nu} z_5^2 z_6+\zeta^{3 \nu} z_6^2 z_4)\\
  &+2 \sqrt{-13-2 \sqrt{13}} (\zeta^{3 \nu} z_1 z_2^2+\zeta^{\nu} z_2 z_3^2+\zeta^{9 \nu} z_3 z_1^2)\\
  &-2 \sqrt{-13+2 \sqrt{13}} (\zeta^{10 \nu} z_4 z_5^2+\zeta^{12 \nu} z_5 z_6^2+\zeta^{4 \nu} z_6 z_4^2)\\
  &+2 \sqrt{\frac{-13-3 \sqrt{13}}{2}} (\zeta^{7 \nu} z_1^2 z_4+\zeta^{11 \nu} z_2^2 z_5+\zeta^{8 \nu} z_3^2 z_6)+\\
  &-2 \sqrt{\frac{-13+3 \sqrt{13}}{2}} (\zeta^{6 \nu} z_1 z_4^2+\zeta^{2 \nu} z_2 z_5^2+\zeta^{5 \nu} z_3 z_6^2)+\\
  &+\sqrt{-13-2 \sqrt{13}} (\zeta^{3 \nu} z_1^2 z_5+\zeta^{\nu} z_2^2 z_6+\zeta^{9 \nu} z_3^2 z_4)+\\
  &+\sqrt{-13+2 \sqrt{13}} (\zeta^{10 \nu} z_2 z_4^2+\zeta^{12 \nu} z_3 z_5^2+\zeta^{4 \nu} z_1 z_6^2)+\\
  &+\sqrt{\frac{-13+3 \sqrt{13}}{2}} (\zeta^{6 \nu} z_1^2 z_6+\zeta^{2 \nu} z_2^2 z_4+\zeta^{5 \nu} z_3^2 z_5)+\\
  &+\sqrt{\frac{-13-3 \sqrt{13}}{2}} (\zeta^{7 \nu} z_3 z_4^2+\zeta^{11 \nu} z_1 z_5^2+\zeta^{8 \nu} z_2 z_6^2)+\\
  &+[2(\theta_1-\theta_3)-3(\theta_2-\theta_4)] z_1 z_2 z_3+[2(\theta_4-\theta_2)-3(\theta_1-\theta_3)] z_4 z_5 z_6+
\endaligned$$
$$\aligned
  &-\sqrt{\frac{-13-3 \sqrt{13}}{2}} (\zeta^{11 \nu} z_1 z_2 z_4+\zeta^{8 \nu} z_2 z_3 z_5+\zeta^{7 \nu} z_1 z_3 z_6)+\\
  &+\sqrt{\frac{-13+3 \sqrt{13}}{2}} (\zeta^{2 \nu} z_1 z_4 z_5+\zeta^{5 \nu} z_2 z_5 z_6+\zeta^{6 \nu} z_3 z_4 z_6)+\\
  &-3 \sqrt{\frac{-13-3 \sqrt{13}}{2}} (\zeta^{7 \nu} z_1 z_2 z_5+\zeta^{11 \nu} z_2 z_3 z_6+\zeta^{8 \nu} z_1 z_3 z_4)+\\
  &+3 \sqrt{\frac{-13+3 \sqrt{13}}{2}} (\zeta^{6 \nu} z_2 z_4 z_5+\zeta^{2 \nu} z_3 z_5 z_6+\zeta^{5 \nu} z_1 z_4 z_6)+\\
  &-\sqrt{-13+2 \sqrt{13}} (\zeta^{10 \nu} z_1 z_2 z_6+\zeta^{4 \nu} z_1 z_3 z_5+\zeta^{12 \nu} z_2 z_3 z_4)+\\
  &+\sqrt{-13-2 \sqrt{13}} (\zeta^{3 \nu} z_3 z_4 z_5+\zeta^{9 \nu} z_2 z_4 z_6+\zeta^{\nu} z_1 z_5 z_6).
\endaligned$$

  This leads us to define the following senary cubic forms (cubic forms in six variables):
$$\left\{\aligned
  {\Bbb D}_0 &=z_1 z_2 z_3,\\
  {\Bbb D}_1 &=2 z_2 z_3^2+z_2^2 z_6-z_4^2 z_5+z_1 z_5 z_6,\\
  {\Bbb D}_2 &=-z_6^3+z_2^2 z_4-2 z_2 z_5^2+z_1 z_4 z_5+3 z_3 z_5 z_6,\\
  {\Bbb D}_3 &=2 z_1 z_2^2+z_1^2 z_5-z_4 z_6^2+z_3 z_4 z_5,\\
  {\Bbb D}_4 &=-z_2^2 z_3+z_1 z_6^2-2 z_4^2 z_6-z_1 z_3 z_5,\\
  {\Bbb D}_5 &=-z_4^3+z_3^2 z_5-2 z_3 z_6^2+z_2 z_5 z_6+3 z_1 z_4 z_6,\\
  {\Bbb D}_6 &=-z_5^3+z_1^2 z_6-2 z_1 z_4^2+z_3 z_4 z_6+3 z_2 z_4 z_5,\\
  {\Bbb D}_7 &=-z_2^3+z_3 z_4^2-z_1 z_3 z_6-3 z_1 z_2 z_5+2 z_1^2 z_4,\\
  {\Bbb D}_8 &=-z_1^3+z_2 z_6^2-z_2 z_3 z_5-3 z_1 z_3 z_4+2 z_3^2 z_6,\\
  {\Bbb D}_9 &=2 z_1^2 z_3+z_3^2 z_4-z_5^2 z_6+z_2 z_4 z_6,\\
  {\Bbb D}_{10} &=-z_1 z_3^2+z_2 z_4^2-2 z_4 z_5^2-z_1 z_2 z_6,\\
  {\Bbb D}_{11} &=-z_3^3+z_1 z_5^2-z_1 z_2 z_4-3 z_2 z_3 z_6+2 z_2^2 z_5,\\
  {\Bbb D}_{12} &=-z_1^2 z_2+z_3 z_5^2-2 z_5 z_6^2-z_2 z_3 z_4,\\
  {\Bbb D}_{\infty}&=z_4 z_5 z_6.
\endaligned\right.\eqno{(4.1)}$$

  Let
$$r_0=2(\theta_1-\theta_3)-3(\theta_2-\theta_4), \quad
  r_{\infty}=2(\theta_4-\theta_2)-3(\theta_1-\theta_3),$$
and
$$r_1=\sqrt{-13-2 \sqrt{13}},
  r_2=\sqrt{\frac{-13+3 \sqrt{13}}{2}},
  r_3=\sqrt{-13+2 \sqrt{13}},
  r_4=\sqrt{\frac{-13-3 \sqrt{13}}{2}}.$$
Then
$$\aligned
  &-13 \sqrt{13} ST^{\nu}({\Bbb D}_0)\\
 =&r_0 {\Bbb D}_0+r_1 \zeta^{\nu} {\Bbb D}_1+r_2 \zeta^{2 \nu} {\Bbb D}_2+
   r_1 \zeta^{3 \nu} {\Bbb D}_3+r_3 \zeta^{4 \nu} {\Bbb D}_4+
   r_2 \zeta^{5 \nu} {\Bbb D}_5+r_2 \zeta^{6 \nu} {\Bbb D}_6+\\
  &+r_4 \zeta^{7 \nu} {\Bbb D}_7+r_4 \zeta^{8 \nu} {\Bbb D}_8
   +r_1 \zeta^{9 \nu} {\Bbb D}_9+r_3 \zeta^{10 \nu} {\Bbb D}_{10}+r_4 \zeta^{11 \nu} {\Bbb D}_{11}+
   r_3 \zeta^{12 \nu} {\Bbb D}_{12}+r_{\infty} {\Bbb D}_{\infty}.
\endaligned$$
$$\aligned
  -13 \sqrt{13} S({\Bbb D}_{\infty})
 =&r_{\infty} {\Bbb D}_0-r_3 {\Bbb D}_1-r_4 {\Bbb D}_2-r_3 {\Bbb D}_3+r_1 {\Bbb D}_4
  -r_4 {\Bbb D}_5-r_4 {\Bbb D}_6+\\
  &+r_2 {\Bbb D}_7+r_2 {\Bbb D}_8-r_3 {\Bbb D}_9+r_1 {\Bbb D}_{10}+r_2 {\Bbb D}_{11}+r_1 {\Bbb D}_{12}
   -r_0 {\Bbb D}_{\infty}.
\endaligned$$
Moreover, we have
$$\aligned
  -13 \sqrt{13} S({\Bbb D}_1)
 =&13 r_1 {\Bbb D}_0+q_1 {\Bbb D}_1+q_2 {\Bbb D}_2+q_3 {\Bbb D}_3+
   q_4 {\Bbb D}_4+q_5 {\Bbb D}_5+q_6 {\Bbb D}_6+\\
  &+q_7 {\Bbb D}_7+q_8 {\Bbb D}_8+q_9 {\Bbb D}_9+q_{10} {\Bbb D}_{10}
   +q_{11} {\Bbb D}_{11}+q_{12} {\Bbb D}_{12}-13 r_3 {\Bbb
   D}_{\infty},
\endaligned$$
where
$$\left\{\aligned
  q_1 &=-2 (\zeta-\zeta^{12})-2 (\zeta^5-\zeta^8)+6 (\zeta^3-\zeta^{10})
        -(\zeta^2-\zeta^{11})+4 (\zeta^9-\zeta^4)+2 (\zeta^6-\zeta^7),\\
  q_2 &=-4 (\zeta-\zeta^{12})+3 (\zeta^5-\zeta^8)+3 (\zeta^3-\zeta^{10})
        -(\zeta^2-\zeta^{11})-2 (\zeta^9-\zeta^4),\\
  q_3 &=6 (\zeta-\zeta^{12})-(\zeta^5-\zeta^8)+4 (\zeta^3-\zeta^{10})+
        2 (\zeta^2-\zeta^{11})-2 (\zeta^9-\zeta^4)-2 (\zeta^6-\zeta^7),\\
  q_4 &=-2 (\zeta-\zeta^{12})+4 (\zeta^5-\zeta^8)+2 (\zeta^3-\zeta^{10})
        -2 (\zeta^2-\zeta^{11})+(\zeta^9-\zeta^4)+6 (\zeta^6-\zeta^7),\\
  q_5 &=-2 (\zeta-\zeta^{12})-4 (\zeta^3-\zeta^{10})+3 (\zeta^2-\zeta^{11})
        +3 (\zeta^9-\zeta^4)-(\zeta^6-\zeta^7),\\
  q_6 &=3 (\zeta-\zeta^{12})-(\zeta^5-\zeta^8)-2 (\zeta^3-\zeta^{10})
        -4 (\zeta^9-\zeta^4)+3 (\zeta^6-\zeta^7),\\
  q_7 &=(\zeta-\zeta^{12})+3 (\zeta^5-\zeta^8)-2 (\zeta^2-\zeta^{11})
        -3 (\zeta^9-\zeta^4)-4 (\zeta^6-\zeta^7),\\
  q_8 &=-2 (\zeta^5-\zeta^8)-3 (\zeta^3-\zeta^{10})-4 (\zeta^2-\zeta^{11})
        +(\zeta^9-\zeta^4)+3 (\zeta^6-\zeta^7),\\
  q_9 &=4 (\zeta-\zeta^{12})+2 (\zeta^5-\zeta^8)-2 (\zeta^3-\zeta^{10})
        -2 (\zeta^2-\zeta^{11})+6 (\zeta^9-\zeta^4)-(\zeta^6-\zeta^7),\\
  q_{10} &=(\zeta-\zeta^{12})+6 (\zeta^5-\zeta^8)-2 (\zeta^3-\zeta^{10})
          +4 (\zeta^2-\zeta^{11})+2 (\zeta^9-\zeta^4)-2 (\zeta^6-\zeta^7),\\
  q_{11} &=-3 (\zeta-\zeta^{12})-4 (\zeta^5-\zeta^8)+(\zeta^3-\zeta^{10})
           +3 (\zeta^2-\zeta^{11})-2 (\zeta^6-\zeta^7),\\
  q_{12} &=2 (\zeta-\zeta^{12})-2 (\zeta^5-\zeta^8)+(\zeta^3-\zeta^{10})
           +6 (\zeta^2-\zeta^{11})-2 (\zeta^9-\zeta^4)+4 (\zeta^6-\zeta^7).
\endaligned\right.\eqno{(4.2)}$$
Similarly, we obtain
$$\aligned
   -13 \sqrt{13} S({\Bbb D}_2)
 =&26 r_2 {\Bbb D}_0+2 q_2 {\Bbb D}_1-q_4 {\Bbb D}_2+2 q_6 {\Bbb D}_3+
   2 q_8 {\Bbb D}_4-q_{10} {\Bbb D}_5-q_{12} {\Bbb D}_6+\\
  &+q_1 {\Bbb D}_7+q_3 {\Bbb D}_8+2 q_5 {\Bbb D}_9+2 q_7 {\Bbb D}_{10}+
   q_9 {\Bbb D}_{11}+2 q_{11} {\Bbb D}_{12}-26 r_4 {\Bbb D}_{\infty}.
\endaligned$$
$$\aligned
  -13 \sqrt{13} S({\Bbb D}_3)
 =&13 r_1 {\Bbb D}_0+q_3 {\Bbb D}_1+q_6 {\Bbb D}_2+q_9 {\Bbb D}_3+
   q_{12} {\Bbb D}_4+q_2 {\Bbb D}_5+q_5 {\Bbb D}_6+\\
  &+q_8 {\Bbb D}_7+q_{11} {\Bbb D}_8+q_1 {\Bbb D}_9+q_4 {\Bbb D}_{10}
   +q_7 {\Bbb D}_{11}+q_{10} {\Bbb D}_{12}-13 r_3 {\Bbb D}_{\infty},
\endaligned$$
$$\aligned
  -13 \sqrt{13} S({\Bbb D}_4)
 =&13 r_3 {\Bbb D}_0+q_4 {\Bbb D}_1+q_8 {\Bbb D}_2+q_{12} {\Bbb D}_3-
   q_3 {\Bbb D}_4+q_7 {\Bbb D}_5+q_{11} {\Bbb D}_6+\\
  &-q_2 {\Bbb D}_7-q_6 {\Bbb D}_8+q_{10} {\Bbb D}_9-q_1 {\Bbb D}_{10}
   -q_5 {\Bbb D}_{11}-q_9 {\Bbb D}_{12}+13 r_1 {\Bbb D}_{\infty},
\endaligned$$
$$\aligned
   -13 \sqrt{13} S({\Bbb D}_5)
 =&26 r_2 {\Bbb D}_0+2 q_5 {\Bbb D}_1-q_{10} {\Bbb D}_2+2 q_2 {\Bbb D}_3+
   2 q_7 {\Bbb D}_4-q_{12} {\Bbb D}_5-q_4 {\Bbb D}_6+\\
  &+q_9 {\Bbb D}_7+q_1 {\Bbb D}_8+2 q_6 {\Bbb D}_9+2 q_{11} {\Bbb D}_{10}+
   q_3 {\Bbb D}_{11}+2 q_8 {\Bbb D}_{12}-26 r_4 {\Bbb D}_{\infty}.
\endaligned$$
$$\aligned
   -13 \sqrt{13} S({\Bbb D}_6)
 =&26 r_2 {\Bbb D}_0+2 q_6 {\Bbb D}_1-q_{12} {\Bbb D}_2+2 q_5 {\Bbb D}_3+
   2 q_{11} {\Bbb D}_4-q_4 {\Bbb D}_5-q_{10} {\Bbb D}_6+\\
  &+q_3 {\Bbb D}_7+q_9 {\Bbb D}_8+2 q_2 {\Bbb D}_9+2 q_8 {\Bbb D}_{10}+
   q_1 {\Bbb D}_{11}+2 q_7 {\Bbb D}_{12}-26 r_4 {\Bbb D}_{\infty}.
\endaligned$$
$$\aligned
   -13 \sqrt{13} S({\Bbb D}_7)
 =&26 r_4 {\Bbb D}_0+2 q_7 {\Bbb D}_1+q_1 {\Bbb D}_2+2 q_8 {\Bbb D}_3-
   2 q_2 {\Bbb D}_4+q_9 {\Bbb D}_5+q_3 {\Bbb D}_6+\\
  &+q_{10} {\Bbb D}_7+q_4 {\Bbb D}_8+2 q_{11} {\Bbb D}_9-2 q_5 {\Bbb D}_{10}+
   q_{12} {\Bbb D}_{11}-2 q_6 {\Bbb D}_{12}+26 r_2 {\Bbb D}_{\infty}.
\endaligned$$
$$\aligned
   -13 \sqrt{13} S({\Bbb D}_8)
 =&26 r_4 {\Bbb D}_0+2 q_8 {\Bbb D}_1+q_3 {\Bbb D}_2+2 q_{11} {\Bbb D}_3-
   2 q_6 {\Bbb D}_4+q_1 {\Bbb D}_5+q_9 {\Bbb D}_6+\\
  &+q_4 {\Bbb D}_7+q_{12} {\Bbb D}_8+2 q_7 {\Bbb D}_9-2 q_2 {\Bbb D}_{10}+
   q_{10} {\Bbb D}_{11}-2 q_5 {\Bbb D}_{12}+26 r_2 {\Bbb D}_{\infty}.
\endaligned$$
$$\aligned
  -13 \sqrt{13} S({\Bbb D}_9)
 =&13 r_1 {\Bbb D}_0+q_9 {\Bbb D}_1+q_5 {\Bbb D}_2+q_1 {\Bbb D}_3+
   q_{10} {\Bbb D}_4+q_6 {\Bbb D}_5+q_2 {\Bbb D}_6+\\
  &+q_1 {\Bbb D}_7+q_7 {\Bbb D}_8+q_8 {\Bbb D}_9+q_{12} {\Bbb D}_{10}
   +q_8 {\Bbb D}_{11}+q_4 {\Bbb D}_{12}-13 r_3 {\Bbb D}_{\infty},
\endaligned$$
$$\aligned
  -13 \sqrt{13} S({\Bbb D}_{10})
 =&13 r_3 {\Bbb D}_0+q_{10} {\Bbb D}_1+q_7 {\Bbb D}_2+q_4 {\Bbb D}_3-
   q_1 {\Bbb D}_4+q_{11} {\Bbb D}_5+q_8 {\Bbb D}_6+\\
  &-q_5 {\Bbb D}_7-q_2 {\Bbb D}_8+q_{12} {\Bbb D}_9-q_9 {\Bbb D}_{10}
   -q_6 {\Bbb D}_{11}-q_3 {\Bbb D}_{12}+13 r_1 {\Bbb D}_{\infty},
\endaligned$$
$$\aligned
   -13 \sqrt{13} S({\Bbb D}_{11})
 =&26 r_4 {\Bbb D}_0+2 q_{11} {\Bbb D}_1+q_9 {\Bbb D}_2+2 q_7 {\Bbb D}_3-
   2 q_5 {\Bbb D}_4+q_3 {\Bbb D}_5+q_1 {\Bbb D}_6+\\
  &+q_{12} {\Bbb D}_7+q_{10} {\Bbb D}_8+2 q_8 {\Bbb D}_9-2 q_6 {\Bbb D}_{10}+
   q_4 {\Bbb D}_{11}-2 q_2 {\Bbb D}_{12}+26 r_2 {\Bbb D}_{\infty}.
\endaligned$$
$$\aligned
  -13 \sqrt{13} S({\Bbb D}_{12})
 =&13 r_3 {\Bbb D}_0+q_{12} {\Bbb D}_1+q_{11} {\Bbb D}_2+q_{10} {\Bbb D}_3-
   q_9 {\Bbb D}_4+q_8 {\Bbb D}_5+q_7 {\Bbb D}_6+\\
  &-q_6 {\Bbb D}_7-q_5 {\Bbb D}_8+q_4 {\Bbb D}_9-q_3 {\Bbb D}_{10}
   -q_2 {\Bbb D}_{11}-q_1 {\Bbb D}_{12}+13 r_1 {\Bbb D}_{\infty}.
\endaligned$$
Hence, we get the element $\widehat{S}$ which is induced from the
action of $S$ on the basis $({\Bbb D}_0, \cdots, {\Bbb
D}_{\infty})$:
$$\widehat{S}=-\frac{1}{13 \sqrt{13}} \left(\matrix S_1 & S_2\\
               S_3 & S_4 \endmatrix\right),\eqno{(4.3)}$$
with
$$S_1=\left(\matrix
  r_0 & r_1 & r_2 & r_1 & r_3 & r_2 & r_2\\
  13 r_1 & q_1 & q_2 & q_3 & q_4 & q_5 & q_6\\
  26 r_2 & 2 q_2 & -q_4 & 2 q_6 & 2 q_8 & -q_{10} & -q_{12}\\
  13 r_1 & q_3 & q_6 & q_9 & q_{12} & q_2 & q_5\\
  13 r_3 & q_4 & q_8 & q_{12} & -q_3 & q_7 & q_{11}\\
  26 r_2 & 2 q_5 & -q_{10} & 2 q_2 & 2 q_7 & -q_{12} & -q_4\\
  26 r_2 & 2 q_6 & -q_{12} & 2 q_5 & 2 q_{11} & -q_4 & -q_{10}
\endmatrix\right),\eqno{(4.4)}$$
$$S_2=\left(\matrix
  r_4 & r_4 & r_1 & r_3 & r_4 & r_3 & r_{\infty}\\
  q_7 & q_8 & q_9 & q_{10} & q_{11} & q_{12} & -13 r_3\\
  q_1 & q_3 & 2 q_5 & 2 q_7 & q_9 & 2 q_{11} & -26 r_4\\
  q_8 & q_{11} & q_1 & q_4 & q_7 & q_{10} & -13 r_3\\
  -q_2 & -q_6 & q_{10} & -q_1 & -q_5 & -q_9 & 13 r_1\\
  q_9 & q_1 & 2 q_6 & 2 q_{11} & q_3 & 2 q_8 & -26 r_4\\
  q_3 & q_9 & 2 q_2 & 2 q_8 & q_1 & 2 q_7 & -26 r_4
\endmatrix\right),\eqno{(4.5)}$$
$$S_3=\left(\matrix
  26 r_4 & 2 q_7 & q_1 & 2 q_8 & -2 q_2 & q_9 & q_3\\
  26 r_4 & 2 q_8 & q_3 & 2 q_{11} & -2 q_6 & q_1 & q_9\\
  13 r_1 & q_9 & q_5 & q_1 & q_{10} & q_6 & q_2\\
  13 r_3 & q_{10} & q_7 & q_4 & -q_1 & q_{11} & q_8\\
  26 r_4 & 2 q_{11} & q_9 & 2 q_7 & -2 q_5 & q_3 & q_1\\
  13 r_3 & q_{12} & q_{11} & q_{10} & -q_9 & q_8 & q_7\\
  r_{\infty} & -r_3 & -r_4 & -r_3 & r_1 & -r_4 & -r_4
\endmatrix\right),\eqno{(4.6)}$$
$$S_4=\left(\matrix
  q_{10} & q_4 & 2 q_{11} & -2 q_5 & q_{12} & -2 q_6 & 26 r_2\\
  q_4 & q_{12} & 2 q_7 & -2 q_2 & q_{10} & -2 q_5 & 26 r_2\\
  q_{11} & q_7 & q_3 & q_{12} & q_8 & q_4 & -13 r_3\\
  -q_5 & -q_2 & q_{12} & -q_9 & -q_6 & -q_3 & 13 r_1\\
  q_{12} & q_{10} & 2 q_8 & -2 q_6 & q_4 & -2 q_2 & 26 r_2\\
  -q_6 & -q_5 & q_4 & -q_3 & -q_2 & -q_1 & 13 r_1\\
  r_2 & r_2 & -r_3 & r_1 & r_2 & r_1 & -r_0
\endmatrix\right).\eqno{(4.7)}$$
Similarly, the element $\widehat{T}$ is induced from the action of
$T$ on the basis $({\Bbb D}_0, \cdots, {\Bbb D}_{\infty})$:
$$\widehat{T}=\text{diag}(1, \zeta, \zeta^2, \zeta^3, \zeta^4,
              \zeta^5, \zeta^6, \zeta^7, \zeta^8, \zeta^9, \zeta^{10},
              \zeta^{11}, \zeta^{12}, 1).\eqno{(4.8)}$$
We have
$$\text{Tr}(\widehat{S})=0, \quad \text{Tr}(\widehat{T})=1, \quad
  \text{Tr}(\widehat{S} \widehat{T})=-2.\eqno{(4.9)}$$
Hence, this fourteen-dimensional representation corresponds to the
character $\chi_{15}$ of $PSL(2, 13)$ in \cite{CC}. It is defined
over the cyclotomic field ${\Bbb Q}(\zeta)$.

  Let
$$\delta_{\infty}(z_1, z_2, z_3, z_4, z_5, z_6)=13^2 (z_1^2 z_2^2 z_3^2+z_4^2 z_5^2 z_6^2)\eqno{(4.10)}$$
and
$$\delta_{\nu}(z_1, z_2, z_3, z_4, z_5, z_6)=\delta_{\infty}(ST^{\nu}(z_1, z_2, z_3, z_4, z_5, z_6))\eqno{(4.11)}$$
for $\nu=0, 1, \cdots, 12$. Then
$$\delta_{\nu}=13^2 ST^{\nu}({\Bbb G}_0)=-13 {\Bbb G}_0+\zeta^{\nu} {\Bbb G}_1
  +\zeta^{2 \nu} {\Bbb G}_2+\cdots+\zeta^{12 \nu} {\Bbb G}_{12},\eqno{(4.12)}$$
where the senary sextic forms (i.e., sextic forms in six
variables) are given as follows:
$$\left\{\aligned
  {\Bbb G}_0 &={\Bbb D}_0^2+{\Bbb D}_{\infty}^2,\\
  {\Bbb G}_1 &=-{\Bbb D}_7^2+2 {\Bbb D}_0 {\Bbb D}_1+10 {\Bbb D}_{\infty} {\Bbb D}_1
               +2 {\Bbb D}_2 {\Bbb D}_{12}-2 {\Bbb D}_3 {\Bbb D}_{11}-4 {\Bbb D}_4 {\Bbb D}_{10}
               -2 {\Bbb D}_9 {\Bbb D}_5,\\
  {\Bbb G}_2 &=-2 {\Bbb D}_1^2-4 {\Bbb D}_0 {\Bbb D}_2+6 {\Bbb D}_{\infty} {\Bbb D}_2
               -2 {\Bbb D}_4 {\Bbb D}_{11}+2 {\Bbb D}_5 {\Bbb D}_{10}-2 {\Bbb D}_6 {\Bbb D}_9
               -2 {\Bbb D}_7 {\Bbb D}_8,\\
  {\Bbb G}_3 &=-{\Bbb D}_8^2+2 {\Bbb D}_0 {\Bbb D}_3+10 {\Bbb D}_{\infty} {\Bbb D}_3
               +2 {\Bbb D}_6 {\Bbb D}_{10}-2 {\Bbb D}_9 {\Bbb D}_7-4 {\Bbb D}_{12} {\Bbb D}_4
               -2 {\Bbb D}_1 {\Bbb D}_2,\\
  {\Bbb G}_4 &=-{\Bbb D}_2^2+10 {\Bbb D}_0 {\Bbb D}_4-2 {\Bbb D}_{\infty} {\Bbb D}_4
               +2 {\Bbb D}_5 {\Bbb D}_{12}-2 {\Bbb D}_9 {\Bbb D}_8-4 {\Bbb D}_1 {\Bbb D}_3
               -2 {\Bbb D}_{10} {\Bbb D}_7,\\
  {\Bbb G}_5 &=-2 {\Bbb D}_9^2-4 {\Bbb D}_0 {\Bbb D}_5+6 {\Bbb D}_{\infty} {\Bbb D}_5
               -2 {\Bbb D}_{10} {\Bbb D}_8+2 {\Bbb D}_6 {\Bbb D}_{12}-2 {\Bbb D}_2 {\Bbb D}_3
               -2 {\Bbb D}_{11} {\Bbb D}_7,\\
  {\Bbb G}_6 &=-2 {\Bbb D}_3^2-4 {\Bbb D}_0 {\Bbb D}_6+6 {\Bbb D}_{\infty} {\Bbb D}_6
               -2 {\Bbb D}_{12} {\Bbb D}_7+2 {\Bbb D}_2 {\Bbb D}_4-2 {\Bbb D}_5 {\Bbb D}_1
               -2 {\Bbb D}_8 {\Bbb D}_{11},\\
  {\Bbb G}_7 &=-2 {\Bbb D}_{10}^2+6 {\Bbb D}_0 {\Bbb D}_7+4 {\Bbb D}_{\infty} {\Bbb D}_7
               -2 {\Bbb D}_1 {\Bbb D}_6-2 {\Bbb D}_2 {\Bbb D}_5-2 {\Bbb D}_8 {\Bbb D}_{12}
               -2 {\Bbb D}_9 {\Bbb D}_{11},\\
  {\Bbb G}_8 &=-2 {\Bbb D}_4^2+6 {\Bbb D}_0 {\Bbb D}_8+4 {\Bbb D}_{\infty} {\Bbb D}_8
               -2 {\Bbb D}_3 {\Bbb D}_5-2 {\Bbb D}_6 {\Bbb D}_2-2 {\Bbb D}_{11} {\Bbb D}_{10}
               -2 {\Bbb D}_1 {\Bbb D}_7,\\
  {\Bbb G}_9 &=-{\Bbb D}_{11}^2+2 {\Bbb D}_0 {\Bbb D}_9+10 {\Bbb D}_{\infty} {\Bbb D}_9
               +2 {\Bbb D}_5 {\Bbb D}_4-2 {\Bbb D}_1 {\Bbb D}_8-4 {\Bbb D}_{10} {\Bbb D}_{12}
               -2 {\Bbb D}_3 {\Bbb D}_6,\\
  {\Bbb G}_{10} &=-{\Bbb D}_5^2+10 {\Bbb D}_0 {\Bbb D}_{10}-2 {\Bbb D}_{\infty} {\Bbb D}_{10}
               +2 {\Bbb D}_6 {\Bbb D}_4-2 {\Bbb D}_3 {\Bbb D}_7-4 {\Bbb D}_9 {\Bbb D}_1
               -2 {\Bbb D}_{12} {\Bbb D}_{11},\\
  {\Bbb G}_{11} &=-2 {\Bbb D}_{12}^2+6 {\Bbb D}_0 {\Bbb D}_{11}+4 {\Bbb D}_{\infty} {\Bbb D}_{11}
               -2 {\Bbb D}_9 {\Bbb D}_2-2 {\Bbb D}_5 {\Bbb D}_6-2 {\Bbb D}_7 {\Bbb D}_4
               -2 {\Bbb D}_3 {\Bbb D}_8,\\
  {\Bbb G}_{12} &=-{\Bbb D}_6^2+10 {\Bbb D}_0 {\Bbb D}_{12}-2 {\Bbb D}_{\infty} {\Bbb D}_{12}
               +2 {\Bbb D}_2 {\Bbb D}_{10}-2 {\Bbb D}_1 {\Bbb D}_{11}-4 {\Bbb D}_3 {\Bbb D}_9
               -2 {\Bbb D}_4 {\Bbb D}_8.
\endaligned\right.\eqno{(4.13)}$$
We have that ${\Bbb G}_0$ is invariant under the action of
$\langle H, T \rangle$, a maximal subgroup of order $78$ of $G$
with index $14$. Note that $\delta_{\infty}$, $\delta_{\nu}$ for
$\nu=0, \cdots, 12$ form an algebraic equation of degree fourteen.
However, we have $\delta_{\infty}+\sum_{\nu=0}^{12}
\delta_{\nu}=0$. Hence, it is not the Jacobian equation of degree
fourteen! We call it exotic modular equation of degree fourteen.
We have
$$\delta_{\infty}^2+\sum_{\nu=0}^{12} \delta_{\nu}^2
 =26 (7 \cdot 13^2 {\Bbb G}_0^2+{\Bbb G}_1 {\Bbb G}_{12}+{\Bbb G}_ 2
  {\Bbb G}_{11}+\cdots+{\Bbb G}_6 {\Bbb G}_7).\eqno{(4.14)}$$
This leads us to define
$$\Phi_{12}(z_1, z_2, z_3, z_4, z_5, z_6)
 :=-\frac{1}{26}(7 \cdot 13^2 {\Bbb G}_0^2+{\Bbb G}_1 {\Bbb G}_{12}
   +{\Bbb G}_2 {\Bbb G}_{11}+\cdots+{\Bbb G}_6 {\Bbb G}_7),\eqno{(4.15)}$$

{\bf Theorem 4.1}. {\it The invariant decomposition formula for
the simple group $PSL(2, 13)$ of order $1092$ is given as follows:
$$\left[\frac{\eta^2(z)}{\eta^2(13z)}\right]^5
 =\frac{\Phi_{12}(a_1(z), a_2(z), a_3(z), a_4(z), a_5(z), a_6(z))}
  {(a_1(z) a_2(z) a_3(z) a_4(z) a_5(z) a_6(z))^2},\eqno{(4.16)}$$
where $\Phi_{12}(z_1, z_2, z_3, z_4, z_5, z_6)$ is an invariant of
degree $12$ associated to $PSL(2, 13)$, and $(z_1 z_2 z_3 z_4 z_5
z_6)^2$ is invariant under the action of the image of a Borel
subgroup of $PSL(2, 13)$ $($i.e. a maximal subgroup of order $78$
of $PSL(2, 13)$$)$.}

{\it Proof}. Let $x_i(z)=\eta(z) a_i(z)$ $(1 \leq i \leq 6)$. Then
$$X(z):=\left(\matrix x_1(z)\\ x_2(z)\\ x_3(z)\\ x_4(z)\\ x_5(z)\\ x_6(z) \endmatrix\right)
       =\eta(z) \left(\matrix a_1(z)\\ a_2(z)\\ a_3(z)\\
        a_4(z)\\ a_5(z)\\ a_6(z) \endmatrix\right)
       =\eta(z) {\bold A}(z).$$
Recall that $\eta(z)$ satisfies the following transformation
formulas
$$\eta(z+1)=e^{\frac{\pi i}{12}} \eta(z), \quad \text{and} \quad
  \eta\left(-\frac{1}{z}\right)=e^{-\frac{\pi i}{4}} \sqrt{z} \eta(z).$$
By Proposition 2.5, we have
$$X(z+1)=e^{-\frac{2 \pi i}{3}} \rho(t) X(z) \quad \text{and} \quad
  X\left(-\frac{1}{z}\right)=z \rho(s) X(z).$$
Hence,
$$X(\gamma(z))=u(\gamma) j(\gamma, z) \rho(\gamma) X(z)
  \quad \text{for $\gamma \in \Gamma(1)$},$$
where $u(\gamma)=1, \omega$ or $\omega^2$ with $\omega=e^{\frac{2
\pi i}{3}}$. Since $\Gamma(13)=\text{ker}$ $\rho$, we have
$$X(\gamma(z))=u(\gamma) j(\gamma, z) X(z) \quad \text{for $\gamma \in \Gamma(13)$}.$$
This means that the functions $x_1(z)$, $x_2(z)$, $x_3(z)$,
$x_4(z)$, $x_5(z)$, $x_6(z)$ are modular forms of weight one for
$\Gamma(13)$ with the same multiplier $u(\gamma)=1, \omega$ or
$\omega^2$. Thus,
$$\Phi_{12}(X(\gamma(z)))=u(\gamma)^{12} j(\gamma, z)^{12} \Phi_{12}(X(z))
 =j(\gamma, z)^{12} \Phi_{12}(X(z)) \quad \text{for $\gamma \in
  \Gamma(13)$}.$$
Moreover, for $\gamma \in \Gamma(1)$,
$$\aligned
  &\Phi_{12}(X(\gamma(z)))=\Phi_{12}(u(\gamma) j(\gamma, z) \rho(\gamma) X(z))\\
 =&u(\gamma)^{12} j(\gamma, z)^{12} \Phi_{12}(\rho(\gamma) X(z))=j(\gamma, z)^{12} \Phi_{12}(\rho(\gamma) X(z)).
\endaligned$$
Note that $\rho(\gamma) \in \langle \rho(s), \rho(t) \rangle=G$
and $\Phi_{12}$ is a $G$-invariant polynomial, we have
$$\Phi_{12}(X(\gamma(z)))=j(\gamma, z)^{12} \Phi_{12}(X(z)), \quad \gamma \in \Gamma(1).$$
This implies that $\Phi_{12}(X(z))$ is a modular form of weight
$12$ for the full modular group $\Gamma(1)$. Moreover, we will
show that it is a cusp form. We have
$$\left\{\aligned
  {\Bbb D}_0 &=q^{\frac{15}{8}} (1+O(q)),\\
  {\Bbb D}_{\infty} &=q^{\frac{7}{8}} (-1+O(q)),\\
  {\Bbb D}_1 &=q^{\frac{99}{104}} (2+O(q)),\\
  {\Bbb D}_2 &=q^{\frac{3}{104}} (-1+O(q)),\\
  {\Bbb D}_3 &=q^{\frac{11}{104}} (1+O(q)),\\
  {\Bbb D}_4 &=q^{\frac{19}{104}} (-2+O(q)),\\
  {\Bbb D}_5 &=q^{\frac{27}{104}} (-1+O(q)),
\endaligned\right. \quad \quad
  \left\{\aligned
  {\Bbb D}_6 &=q^{\frac{35}{104}} (-1+O(q)),\\
  {\Bbb D}_7 &=q^{\frac{43}{104}} (1+O(q)),\\
  {\Bbb D}_8 &=q^{\frac{51}{104}} (3+O(q)),\\
  {\Bbb D}_9 &=q^{\frac{59}{104}} (-2+O(q)),\\
  {\Bbb D}_{10} &=q^{\frac{67}{104}} (1+O(q)),\\
  {\Bbb D}_{11} &=q^{\frac{75}{104}} (-4+O(q)),\\
  {\Bbb D}_{12} &=q^{\frac{83}{104}} (-1+O(q)).
\endaligned\right.$$
Hence,
$$\left\{\aligned
  {\Bbb G}_0 &=q^{\frac{7}{4}} (1+O(q)),\\
  {\Bbb G}_1 &=q^{\frac{86}{104}} (13+O(q)),\\
  {\Bbb G}_2 &=q^{\frac{94}{104}} (-22+O(q)),\\
  {\Bbb G}_3 &=q^{\frac{102}{104}} (-21+O(q)),\\
  {\Bbb G}_4 &=q^{\frac{6}{104}} (-1+O(q)),\\
  {\Bbb G}_5 &=q^{\frac{14}{104}} (2+O(q)),\\
  {\Bbb G}_6 &=q^{\frac{22}{104}} (2+O(q)),
\endaligned\right. \quad \quad
  \left\{\aligned
  {\Bbb G}_7 &=q^{\frac{30}{104}} (-2+O(q)),\\
  {\Bbb G}_8 &=q^{\frac{38}{104}} (-8+O(q)),\\
  {\Bbb G}_9 &=q^{\frac{46}{104}} (6+O(q)),\\
  {\Bbb G}_{10} &=q^{\frac{54}{104}} (1+O(q)),\\
  {\Bbb G}_{11} &=q^{\frac{62}{104}} (-8+O(q)),\\
  {\Bbb G}_{12} &=q^{\frac{70}{104}} (17+O(q)).
\endaligned\right.$$
Therefore,
$$\aligned
  &7 \cdot 13^2 {\Bbb G}_0^2+{\Bbb G}_1 {\Bbb G}_{12}+{\Bbb G}_2
   {\Bbb G}_{11}+\cdots+{\Bbb G}_6 {\Bbb G}_7\\
 =&7 \cdot 13^2 q^{\frac{7}{2}} (1+O(q))+q^{\frac{3}{2}}
   (13 \cdot 17+O(q))+q^{\frac{3}{2}} (22 \cdot 8+O(q))
   +q^{\frac{3}{2}} (-21+O(q))+\\
  &+q^{\frac{1}{2}}(-6+O(q))+q^{\frac{1}{2}}(-16+O(q))+q^{\frac{1}{2}}(-4+O(q))\\
 =&q^{\frac{1}{2}}(-26+O(q)).
\endaligned$$
We have
$$\aligned
 &\Phi_{12}(x_1(z), x_2(z), x_3(z), x_4(z), x_5(z), x_6(z))\\
=&\eta(z)^{12} q^{\frac{1}{2}} (1+O(q))\\
=&q^{\frac{1}{2}} \prod_{n=1}^{\infty} (1-q^n)^{12} \cdot
  q^{\frac{1}{2}} (1+O(q))\\
=&q(1+O(q))
\endaligned$$
is a cusp form of weight $12$ for the full modular group
$\Gamma(1)$. Because every $\Gamma(1)$ cusp form of weight $12$ is
a multiple of $\Delta(z)$, checking the $q^1$ coefficient, we find
that $\Phi_{12}(x_1(z), x_2(z), x_3(z), x_4(z), x_5(z),
x_6(z))=\Delta(z)$. On the other hand, we have
$$x_1(z) x_2(z) x_3(z) x_4(z) x_5(z) x_6(z)=-\eta(z)^7 \eta(13z)^5.$$
This completes the proof of Theorem 4.1.

\flushpar $\qquad \qquad \qquad \qquad \qquad \qquad \qquad \qquad
\qquad \qquad \qquad \qquad \qquad \qquad \qquad \qquad \qquad
\qquad \quad \boxed{}$

  Note that (2.20), (2.21) and (2.22) are equivalent to the following:
$$\aligned
 &a_1(z)^2 a_4(z)^2 a_2(z) a_5(z)+a_2(z)^2 a_5(z)^2 a_3(z) a_6(z)
 +a_3(z)^2 a_6(z)^2 a_1(z) a_4(z)+\\
 &-a_1(z) a_4(z) a_2(z) a_5(z) a_3(z) a_6(z)=-\eta(z)^3 \eta(13z)^3,
\endaligned\eqno{(4.17)}$$
$$\aligned
 &a_1(z)^2 a_4(z)^2 a_3(z) a_6(z)+a_2(z)^2 a_5(z)^2 a_1(z) a_4(z)
 +a_3(z)^2 a_6(z)^2 a_2(z) a_5(z)+\\
 &+4 a_1(z) a_4(z) a_2(z) a_5(z) a_3(z) a_6(z)=\eta(z)^3 \eta(13z)^3,
\endaligned\eqno{(4.18)}$$
$$\aligned
 &a_4(z)^2 a_5(z)^2 a_6(z)^2-a_1(z)^2 a_2(z)^2 a_3(z)^2+3 a_1(z) a_4(z) a_2(z) a_5(z) a_3(z) a_6(z)\\
=&\eta(z)^3 \eta(13z)^3.
\endaligned\eqno{(4.19)}$$
The corresponding polynomials are given by
$$f_6(z_1, z_2, z_3, z_4, z_5, z_6):=z_1^2 z_4^2 z_2 z_5+z_2^2
  z_5^2 z_3 z_6+z_3^2 z_6^2 z_1 z_4-z_1 z_2 z_3 z_4 z_5 z_6,$$
$$g_6(z_1, z_2, z_3, z_4, z_5, z_6):=z_1^2 z_4^2 z_3 z_6+z_2^2
  z_5^2 z_1 z_4+z_3^2 z_6^2 z_2 z_5+4 z_1 z_2 z_3 z_4 z_5 z_6,$$
$$h_6(z_1, z_2, z_3, z_4, z_5, z_6):=z_4^2 z_5^2 z_6^2-z_1^2 z_2^2 z_3^2
  +3 z_1 z_2 z_3 z_4 z_5 z_6.$$
Similar as the above argument in the proof of Theorem 4.1, we have
$$\aligned
 &f_6(X(\gamma(z)))=f_6(u(\gamma) j(\gamma, z) \rho(\gamma) X(z))\\
=&u(\gamma)^6 j(\gamma, z)^6 f_6(\rho(\gamma) X(z))=j(\gamma, z)^6
  f_6(\rho(\gamma) X(z)) \quad \text{for $\gamma \in \Gamma(1)$}.
\endaligned$$
If $f_6$ is a $G$-invariant polynomial, we have
$$f_6(X(\gamma(z)))=j(\gamma, z)^6 f_6(X(z)) \quad \text{for $\gamma \in \Gamma(1)$}.$$
This implies that $f_6(X(z))$ is a modular form of weight $6$ for
the full modular group $\Gamma(1)$. On the other hand, by (4.17),
we have
$$f_6(X(z))=-\eta(z)^9 \eta(13 z)^3=-q^2 \prod_{n=1}^{\infty} (1-q^n)^9
            \prod_{n=1}^{\infty} (1-q^{13n})^3.$$
This shows that $f_6(X(z))$ is a cusp form of weight $6$ for the
full modular group $\Gamma(1)$, but the only such form is zero.
This leads to a contradiction! Therefore, $f_6$ is not a
$G$-invariant polynomial. Similarly, we can prove that $g_6$ and
$h_6$ are not $G$-invariant polynomials.

\vskip 0.5 cm

\centerline{\bf Appendix}

\vskip 0.5 cm

  In this appendix, we give some importance of (1.1) and (1.2) from
another point of view. According to \cite{G2}, we have the
following tables:
$$\text{Table $1$. $\eta$-products in $S_{\frac{1}{2}(p-1)}(p, \chi_p)$ for small $p$}$$
$$\matrix
   p & \text{$\eta$-product cusp forms}\\
   5 & -\\
   7 & \eta(z)^3 \eta(7z)^3\\
  11 & -\\
  13 & \eta(z)^{2n+1} \eta(13z)^{11-2n} \quad (0 \leq n \leq 5)
\endmatrix$$
$$\text{Table $2$. Dimension of certain spaces of cusp forms}$$
$$\matrix
   p & \dim S_{\frac{1}{2}(p-1)}(p, \chi_p)\\
   5 & 0\\
   7 & 1\\
  11 & 3\\
  13 & 6
\endmatrix$$
Note that
$$\aligned
  q^{\frac{19}{24}} \sum_{n=0}^{\infty} p(5n+4) q^n
 &=\frac{1}{\eta(z)^5} \cdot 0+5 \frac{\eta(5z)^5}{\eta(z)^6},\\
  q^{\frac{17}{24}} \sum_{n=0}^{\infty} p(7n+5) q^n
 &=\frac{1}{\eta(z)^7} \cdot 7 \eta(z)^3 \eta(7z)^3+49 \frac{\eta(7z)^7}{\eta(z)^8},
\endaligned\eqno{(A.1)}$$
$$\aligned
  q^{\frac{11}{24}} \sum_{n=0}^{\infty} p(13n+6) q^n
=&\frac{1}{\eta(z)^{13}} [11 \eta(13z) \eta(z)^{11}+36 \cdot
  13 \eta(13z)^3 \eta(z)^9+\\
 &+38 \cdot 13^2 \eta(13z)^5 \eta(z)^7+20 \cdot 13^3 \eta(13z)^7 \eta(z)^5+\\
 &+6 \cdot 13^4 \eta(13z)^9 \eta(z)^3+13^5 \eta(13z)^{11} \eta(z)]+13^5 \frac{\eta(13z)^{13}}{\eta(z)^{14}}.
\endaligned\eqno{(A.2)}$$
We observe that what appear in the right-hand side of (A.1) and
(A.2) are just the base vectors of $S_{\frac{1}{2}(p-1)}(p,
\chi_p)$ for $p=5$, $7$ and $13$.

  In the end, let us recall some basic fact about the Monster finite
simple group ${\Bbb M}$ (see \cite{CN})and its relation with its
type $5B$, $7B$ and $13B$ elements. For each prime $p$ with
$(p-1)|24$ there is a conjugacy class of elements of ${\Bbb M}$,
with centraliser of form $p^{1+2d}.G_p$, where $p.G_p$ is the
centraliser of a corresponding automorphism of the Leech Lattice.
The symbol $p^{1+2d}$ denotes an extraspecial $p$-group, and $2d=
24/(p-1)$. Some maximal $p$-local subgroups of ${\Bbb M}$ are
given as follows (see \cite{CC}):
$$\matrix
  p & \text{structure} & \text{specification}\\
  2 & 2_{+}^{1+24} \cdot Co_1 & N(2B)\\
  3 & 3_{+}^{1+12} \cdot 2 Suz:2 & N(3B)\\
  5 & 5_{+}^{1+6}: 4J_2 \cdot 2 & N(5B)\\
  7 & 7_{+}^{1+4}:(3 \times 2S_7) & N(7B)\\
 13 & 13_{+}^{1+2}:(3 \times 4 S_4) & N(13B)
\endmatrix$$
and
$$\matrix
  p & \text{structure} & \text{specification}\\
  2 & 2^2.2^{11}.2^{22}.(S_3 \times M_{24}) & N(2B^2)\\
  3 & 3^2.3^5.3^{10}.(M_{11} \times 2S_4) & N(3B^2)\\
  5 & 5^2.5^2.5^4.(S_3 \times GL_2(5)) & N(5B^2)\\
  7 & 7^2.7.7^2:GL_2(7) & N(7B^2)\\
 13 & 13^2:4 L_2(13) \cdot 2 & N(13B^2)
\endmatrix$$

{\smc Department of Mathematics, Peking University}

{\smc Beijing 100871, P. R. China}

{\it E-mail address}: yanglei\@math.pku.edu.cn
\vskip 1.5 cm
\Refs

\item{[AB]} {\smc S. Ahlgren and M. Boylan}, Arithmetic properties
            of the partition function, Invent. Math. {\bf 153} (2003),
            487-502.

\item{[AG]} {\smc G. E. Andrews and F. G. Garvan}, Dyson's crank
             of a partition, Bull. Amer. Math. Soc. {\bf 18}
             (1988), 167-171.

\item{[AB1]} {\smc G. E. Andrews and B. C. Berndt}, {\it
             Ramanujan's Lost Notebook}, Part I, Springer, 2005.

\item{[ASD]} {\smc A. O. L. Atkin and H. P. F. Swinnerton-Dyer},
             Some properties of partitions, Proc. London Math. Soc.
             {\bf 66} (1954), 84-106.

\item{[B1]} {\smc B. C. Berndt}, {\it Ramanujan's Notebooks}, Part
             III, Springer-Verlag, 1991.

\item{[B2]} {\smc B. C. Berndt}, {\it Ramanujan's Notebooks}, Part
             IV, Springer-Verlag, 1994.

\item{[BCCL]} {\smc B. C. Berndt, H. H. Chan, S. H. Chan and W.-C.
              Liaw}, Cranks--really, the final problem, Ramanujan J.
              {\bf 23} (2010), 3-15.

\item{[BO]} {\smc B. C. Berndt and K. Ono}, Ramanujan's
             unpublished manuscript on the partition and tau functions with
             proofs and commentary, in: {\it The Andrews
             Festschrift, Seventeen Papers on Classical Number
             Theory and Combinatorics}, 39-110, Springer, 2001.

\item{[Bor]} {\smc R. E. Borcherds}, Monstrous moonshine and
             monstrous Lie superalgebras, Invent. Math. {\bf 109}
             (1992), 405-444.

\item{[BrO]} {\smc K. Bringmann and K. Ono}, Dyson's ranks and
             Maass forms, Ann. of Math. (2) {\bf 171} (2010), 419-449.

\item{[CY]} {\smc I. Chen and N. Yui}, Singular values of Thompson
             series, In: {\it Groups, difference sets, and the Monster
             (Columbus, OH, 1993)}, 255-326, Ohio State Univ. Math. Res.
             Inst. Publ., 4, de Gruyter, Berlin, 1996.

\item{[CN]} {\smc J. Conway and S. Norton}, Monstrous moonshine,
            Bull. London Math. Soc. {\bf 11} (1979), 308-339.

\item{[CC]} {\smc J. H. Conway, R. T. Curtis, S. P. Norton, R.
             A. Parker and R. A. Wilson}, {\it Atlas of Finite Groups,
             Maximal Subgroups and Ordinary Characters for Simple Groups},
             Clarendon Press, Oxford, 1985.

\item{[De1]} {\smc P. Deligne}, Letter to the author, June 29,
             2014. Private communication.

\item{[De2]} {\smc P. Deligne}, Letter to the author, August 7,
             2014. Private communication.

\item{[Du]} {\smc W. Duke}, Continued fractions and modular
             functions, Bull. Amer. Math. Soc. (N.S.) {\bf 42} (2005), 137-162.

\item{[D]} {\smc F. J. Dyson}, Some guesses in the theory of
             partitions, Eureka (Cambridge) {\bf 8} (1944), 10-15.

\item{[El]} {\smc N. D. Elkies}, The Klein quartic in number
             theory, {\it The eightfold way}, 51-101, Math. Sci. Res. Inst.
             Publ., {\bf 35}, Cambridge Univ. Press, Cambridge, 1999.

\item{[Ev]} {\smc R. J. Evans}, Theta function identities, J.
             Math. Anal. Appl. {\bf 147} (1990), 97-121.

\item{[FK]} {\smc H. M. Farkas and I. Kra}, {\it Theta Constants,
            Riemann Surfaces and the Modular Group, An Introduction
            with Applications to Uniformization Theorems, Partition
            Identities and Combinatorial Number Theory}, Graduate
            Studies in Mathematics, {\bf 37}, American Mathematical
            Society, Providence, RI, 2001.

\item{[G1]} {\smc F. G. Garvan}, New combinatorial interpretations
            of Ramanujan's partition congruences mod $5$, $7$ and
            $11$, Trans. Amer. Math. Soc. {\bf 305} (1988),
            47-77.

\item{[G2]} {\smc F. G. Garvan}, Some congruences for partitions
             that are $p$-cores, Proc. London Math. Soc. {\bf 66}
             (1993), 449-478.

\item{[GKS]} {\smc F. Garvan, D. Kim and D. Stanton}, Cranks and
              $t$-cores, Invent. Math. {\bf 101} (1990), 1-17.

\item{[Gr]} {\smc B. H. Gross}, Heegner points and the modular
              curve of prime level, J. Math. Soc. Japan {\bf 39} (1987), 345-362.

\item{[Hu]} {\smc A. Hurwitz}, Grundlagen einer independenten
            Theorie der elliptischen Modulfunktionen und Theorie
            der Multiplikator-Gleichungen erster Stufe, Math. Ann.
            {\bf 18} (1881), 528-592, in: {\it Mathematische
            Werke}, Bd. 1, 1-66, E. Birkh\"{a}user, 1932-33.

\item{[K]} {\smc F. Klein}, {\it Lectures on the Icosahedron and
             the Solution of Equations of the Fifth Degree},
             Translated by G. G. Morrice, second and revised edition,
             Dover Publications, Inc., 1956.

\item{[K1]} {\smc F. Klein}, Ueber die Transformation der
            elliptischen Functionen und die Aufl\"{o}sung der
            Gleichungen f\"{u}nften Grades, Math. Ann. {\bf 14}
            (1879), 111-172, in: {\it Gesammelte Mathematische
            Abhandlungen}, Bd. III, 13-75, Springer-Verlag, Berlin,
            1923.

\item{[K2]} {\smc F. Klein}, Ueber die Transformation siebenter
            Ordnung der elliptischen Functionen, Math. Ann. {\bf
            14} (1879), 428-471, in: {\it Gesammelte Mathematische
            Abhandlungen}, Bd. III, 90-136, Springer-Verlag, Berlin,
            1923.

\item{[K3]} {\smc F. Klein}, Ueber die Aufl\"{o}sung gewisser
            Gleichungen vom siebenten und achten Grade, Math. Ann.
            {\bf 15} (1879), 251-282, in: {\it Gesammelte Mathematische
            Abhandlungen}, Bd. II, 390-438, Springer-Verlag, Berlin,
            1922.

\item{[K4]} {\smc F. Klein}, Ueber die Transformation elfter
            Ordnung der elliptischen Functionen, Math. Ann. {\bf
            15} (1879), 533-555, in: {\it Gesammelte Mathematische
            Abhandlungen}, Bd. III, 140-168, Springer-Verlag, Berlin,
            1923.

\item{[KF1]} {\smc F. Klein and R. Fricke}, {\it Vorlesungen
            \"{u}ber die Theorie der Elliptischen
            Modulfunctionen}, Vol. I, Leipzig, 1890.

\item{[KF2]} {\smc F. Klein and R. Fricke}, {\it Vorlesungen
            \"{u}ber die Theorie der Elliptischen
            Modulfunctionen}, Vol. II, Leipzig, 1892.

\item{[KL]} {\smc D. S. Kubert and S. Lang}, {\it Modular Units},
            Grundlehren der Mathematischen Wissenschaften {\bf
            244}, Springer-Verlag, 1981.

\item{[La]} {\smc G. Lachaud}, Ramanujan modular forms and the
            Klein quartic, Mosc. Math. J. {\bf 5} (2005), 829-856.

\item{[Ma]} {\smc K. Mahlburg}, Partition congruences and the
            Andrews-Garvan-Dyson crank, Proc. Natl. Acad. Sci.
            USA {\bf 102} (2005), 15373-15376.

\item{[O]} {\smc K. Ono}, Distribution of the partition function
            modulo $m$, Ann. of Math. (2) {\bf 151} (2000), 293-307.

\item{[Ra]} {\smc H. Rademacher}, The Ramanujan identities under
            modular substitutions, Trans. Amer. Math. Soc. {\bf 51}
            (1942), 609-636.

\item{[R1]} {\smc S. Ramanujan}, Some properties of $p(n)$, the
            number of partitions of $n$, Proc. Cambridge Phil.
            Soc. {\bf 19} (1919), 207-210, in: {\it Collected
            Papers of Srinivasa Ramanujan}, 210-213, Cambridge,
            1927.

\item{[R2]} {\smc S. Ramanujan}, Congruence properties of
            partitions, Proc. London Math. Soc. {\bf 18} (1920), 19,
            in: {\it Collected Papers of Srinivasa Ramanujan}, 230,
            Cambridge, 1927.

\item{[R3]} {\smc S. Ramanujan}, Congruence properties of
            partitions, Math. Z. {\bf 9} (1921), 147-153, in:
            {\it Collected Papers of Srinivasa Ramanujan}, 232-238,
            Cambridge, 1927.

\item{[R4]} {\smc S. Ramanujan}, Algebraic relations between
            certain infinite products, Proc. London Math. Soc. {\bf 18}
            (1920), in: {\it Collected Papers of Srinivasa Ramanujan}, 231,
            Cambridge, 1927.

\item{[R5]} {\smc S. Ramanujan}, {\it Notebooks of Srinivasa
            Ramanujan}, Vol. I, Springer-Verlag, 1984.

\item{[R6]} {\smc S. Ramanujan}, {\it Notebooks of Srinivasa
            Ramanujan}, Vol. II, Springer-Verlag, 1984.

\item{[R]} {\smc S. Ramanujan}, {\it The Lost Notebook and Other
            Unpublished Papers}, Narosa, New Delhi, 1988.

\item{[Se]} {\smc J.-P. Serre}, {\it A Course in Arithmetic},
           Graduate Texts in Math. {\bf 7}, Springer, New York, 1973.

\item{[S]} {\smc A. Sinkov}, Necessary and sufficient conditions
            for generating certain simple groups by two operators
            of periods two and three, Amer. J. Math. {\bf 59}
            (1937), 67-76.

\item{[W]} {\smc A. Weil}, Sur certains groupes d'op\'{e}rateurs
           unitaires, Acta Math. {\bf 111} (1964), 143-211.

\item{[Y]} {\smc L. Yang}, Exotic arithmetic structure on the
           first Hurwitz triplet. arXiv:1209.1783v5 [math.NT],
           2013.

\item{[Z]} {\smc H. S. Zuckerman}, Identities analogous to
            Ramanujan's identities involving the partition
            function, Duke Math. J. {\bf 5} (1939), 88-110.

\endRefs
\end{document}